\documentclass{amsproc}
\usepackage[latin1]{inputenc}
\usepackage[T1]{fontenc}
\usepackage{mathrsfs,graphics,color}
\usepackage{hyperref}
\def\cprime{$'$} 
\newcommand{\mytypeout}[1]{}
\def\=mod#1[#2]{\equiv#1\,[#2]}
\newcommand{\PT}{{\rm (T)}}
\newcommand{\fset}[1]{\{0,\ldots,#1-1\}}
\newcommand{\lcm}{\mrm{lcm}\,}

\newcommand{\mrm}[1]{\text{\rm #1}}
\newcommand{\swei}[1]{\|#1\|}
\newcommand{\cR}{\mathcal{R}} 
\newcommand{\sB}{\mathscr{B}}
\newcommand{\sK}{\mathscr{K}}
\newcommand{\sT}{\mathscr{T}}
\newcommand{\sE}{\mathscr{E}}
\newcommand{\El}{\mathscr{E}_{\lambda}}
\newcommand{\sM}{\mathscr{M}}

\newcommand{\sMin}{\sM^{m}}
\newcommand{\set}[2]{\{#1\mid\,#2\}} 
\newcommand{\R}{\mathbb{R}}
\newcommand{\N}{\mathbb{N}}

\newcommand{\rmax}{\R_{\max}}
\newcommand{\rmaxb}{\overline{\R}_{\max}}
\newcommand{\new}[1]{{\em #1}\index{#1}}
\newcommand{\tr}{\operatorname{tr}}
\newcommand{\pl}{\flat}
\newcommand{\cyc}[1]{\gamma (#1)}
\newcommand{\cycc}[1]{\sigma (#1)}

\DeclareMathAlphabet{\mathbbold}{U}{bbold}{m}{n}
\newcommand{\zero}{\mathbbold{0}}
\newcommand{\unit}{\mathbbold{1}}
\newcommand{\Al}{A_{\lambda}}
\newtheorem{prop}{Proposition}[section]
\newtheorem{corollary}[prop]{Corollary}
\newtheorem{lemma}[prop]{Lemma}
\newtheorem{assumption}[prop]{Assumption}
\newtheorem{obs}[prop]{Observation}
\newtheorem{theorem}[prop]{Theorem}
\theoremstyle{definition}
\newtheorem{definition}[prop]{Definition}
\theoremstyle{remark}
\newtheorem{example}[prop]{Example}
\newtheorem{remark}[prop]{Remark}
\begin{document}
\title{Discrete max-plus spectral theory}
\author{Marianne Akian}
\address{Marianne Akian: INRIA, Domaine de Voluceau, 
78153 Le Chesnay C\'edex, France}
\email{marianne.akian@inria.fr}
\author{St\'ephane Gaubert}
\address{St\'ephane Gaubert: INRIA, Domaine de Voluceau, 
78153 Le Chesnay C\'edex, France}
\email{stephane.gaubert@inria.fr}
\author{Cormac Walsh}
\address{Cormac Walsh: INRIA, Domaine de Voluceau, 
78153 Le Chesnay C\'edex, France}
\email{cormac.walsh@inria.fr}
\date{October 24, 2003. Revised April 14, 2004.}
\thanks{This work was partially supported by the Erwin Schr\"odinger International Institute for Mathematical Physics (ESI)}
\thanks{This work was done during a post-doctoral stay of the third author
at INRIA, supported by an ERCIM-INRIA fellowship}
\begin{abstract}
We develop a max-plus spectral theory for infinite matrices.
We introduce recurrence and tightness conditions,
under which many results of the finite dimensional
theory, concerning the representation
of eigenvectors and the asymptotic
behavior of powers of matrices, carry over.
We also announce more general representation
results for eigenvectors, which are
obtained by introducing a max-plus analogue
of the Martin boundary.
\end{abstract}
\maketitle
\section{Introduction}
Given a set $S$ and a map $A:S\times S\to \R\cup\{-\infty\},\;
(i,j)\mapsto A_{ij}$,
we consider the following spectral problem:
find $\lambda \in \R\cup\{-\infty\}$ and
$u:S\to\R\cup\{-\infty\},\; i\mapsto u_i$, not identically $-\infty$,
such that
\begin{align*}
\lambda+u_i = \sup_{j\in S} (A_{ij} + u_j)
\qquad\mbox{for all $i\in S$.}
\end{align*}
We say that $\lambda$ is an \new{eigenvalue} of $A$
and that $u$ is an \new{eigenvector} (or a \new{$\lambda$-eigenvector}).
This terminology can be justified by introducing
the max-plus semiring, which is the set $\rmax:=\R\cup\{-\infty\}$
equipped with the addition operation $a\oplus b:=\max(a,b)$
and multiplication operation $a\otimes b:=a+b$.
We write $\zero:=-\infty$ and $\unit:=0$ for the zero and unit
elements of $\rmax$, respectively, and also
denote by $\zero$ the identically $\zero$ vector.

In max-plus notation, the spectral problem
becomes:
\begin{align*}
\text{find $u\in \rmax^S\setminus\{\zero\}$ such that $\lambda u = Au$,}
\end{align*}
where $A:\rmax^S\to (\rmax\cup\{+\infty\})^S$ is the operator defined by
\begin{align}
(Au)_i:=\bigoplus_{j\in S} A_{ij} u_j \enspace .
\label{e-kernel}
\end{align}
Here and in the rest of the paper,
we adopt the usual algebraic conventions, writing for instance $ab$ instead
of $a\otimes b$. We also use the notation $\bigoplus$
for the supremum of an arbitrary family.
 
The operator $A$ is \new{max-plus linear}, meaning that
$A(\alpha u\oplus \beta v)=\alpha Au
\oplus \beta Av$, for all $u,v\in\rmax^S$
and $\alpha,\beta\in \rmax$.
A max-plus ``Riesz representation theorem'',
due to Kolokoltsov and Maslov~\cite{KOLO88,KOLO90,maslovkololtsov95} 
and Akian~\cite{DENSITE} states that under
fairly general conditions, max-plus linear maps
can be represented in the form given in~\eqref{e-kernel} 
(see also~\cite{maxplus91a},\cite{shubin},\cite{kol92},
\cite[Th.~6.5]{bcoq}, \cite[Ch.~8]{singer97}, \cite{agk04}).
The map $(i,j)\mapsto A_{ij}$ is called
the \new{kernel} or \new{matrix} of the operator $A$.
\mytypeout{PUHALSKI}

Max-plus linear operators with kernels
arise as dynamic programming operators
associated to optimal control problems
(see for instance~\cite{romanovski}).
Here, the eigenvectors are the stationary solutions
of the dynamic programming equations and $\lambda$
is the maximal ergodic reward per time unit.
The spectral problem also
arises in the study of discrete event systems~\cite{cohen85a,bcoq},
in statistical mechanics~\cite{griffiths}, in perturbation
problems for eigenvalues and eigenvectors~\cite{ABG96,abg04a},
and in the study of delay systems~\cite{malletparetnussbaum}.

The spectral problem has been much studied
in the ``matrix case'', that is when $S$ is a finite set.
A basic result, which was obtained independently by several
authors, including Cuninghame-Green~\cite{cuning79},
Romanovski\u\i~\cite{romanovski}, and Vorob{\cprime}ev~\cite{vorobyev67},
is that when $A$ is irreducible, it has a unique eigenvalue
which coincides with the maximal circuit mean of $A$.
Gondran and Minoux~\cite{gondran77} and 
Cuninghame-Green~\cite[Th.~24.9]{cuning79} both obtained
a generating family of the eigenspace.
Cohen, Dubois, Quadrat, and Viot~\cite{cohen83} 
showed that the asymptotic behavior of the powers of $A$ can be described
explicitly in terms of the eigenvectors. A useful tool in  this theory
is the \new{critical graph},
which consists of the circuits with maximal circuit mean.
It allows one to determine the dimension
of the eigenspace and the ultimate period of the powers of $A$.
Other references
on the subject are~\cite{Zimmermann.U,gondran84,kim,wende,nussbaum91,bcoq,dudnikov92b,cuning95b,bapat95,maxplus97,bapat98,gondran02}.

The case when $S$ is infinite has also received some
attention~\cite{dudnikov,kol92,lesin,yakovenko,maslovkololtsov95,malletparetnussbaum},
particularly the case when $S$ is compact, $A$ 
leaves the space of continuous functions from $S$ to $\rmax$ invariant,
and the eigenvector is required to be in this space.
For instance, in~\cite[Th.~2]{kol92}, Kolokoltsov characterized the case
where $A$ is a compact map.
In~\cite[Section 2.3]{maslovkololtsov95},
Kolokoltsov and Maslov gave,
under some additional assumptions,
existence results for the eigenvector,
using a Krein-Rutman type fixed point approach.
They also described the eigenspace and gave convergence results
in special cases. Some of the results of these two references
apply more generally to the case where $S$ is locally compact,
provided one considers only eigenvectors that tend
to $\zero$ at infinity.
In~\cite{malletparetnussbaum}, Mallet-Paret and Nussbaum
showed the existence of eigenvectors in the case 
when $S$ is a compact interval of $\R$ and $A$ is a particular 
non-compact map, using measures of non-compactness.

We should also remark that
the continuous time version of the max-plus spectral problem
appears in the study of the Lagrange problem of calculus of variations
and of Hamilton-Jacobi equations, see~\cite{maslov92},~\cite[Ch.~3]{maslovkololtsov95}, and in the related subject of weak
KAM-theory~\cite{mather,mane,fathi97b,fathi03}.

In this paper, we assume that $S$
is an infinite discrete topological space, and 
develop a general max-plus spectral theory in this context.

Section~\ref{sec-graph} is devoted to preliminaries.
In particular, we recall the definition and properties of 
the max-plus analogue of the \new{potential kernel},
$A^*:=I\oplus A\oplus A^2\oplus \cdots$, and of the \new{maximal circuit mean},
$\rho(A)$.
A useful observation, made in~\cite{dudnikov}, is
that every eigenvalue associated to an eigenvector with full
support must be greater than or equal to $\rho(A)$.
Unlike in the finite dimensional case, the spectrum of $A$ may
differ from $\{\rho(A)\}$ even when $A$ is irreducible.

In Sections~\ref{sec-rec} to~\ref{sec-eigstar},
we study the eigenspace associated to the value $\rho(A)$,
extending some results from the finite dimensional theory.

A new feature is the notion of {\em recurrent} nodes and classes,
which replaces the notion of {\em critical}  nodes and classes
appearing in the finite-dimensional theory.
A node is recurrent (Definition~\ref{defi-rec})
if one can return to it with ``normalized'' reward $\unit$.
In algebraic terms, assuming $\rho(A)\in\R$
and setting $\tilde{A}=\rho(A)^{-1} A$ and
$\tilde{A}^+:=\tilde{A}\tilde{A}^*$,
we say that a node $i$ is recurrent if $\tilde{A}^+_{ii}=\unit$.
Associated to any recurrent node $i$ is an eigenvector of $A$,
obtained by taking the $i$th column
of $\tilde{A}^+$ (Proposition~\ref{l2}). 

We then introduce a tightness condition (Property \PT, 
see Definition~\ref{def-ptyt}) which may be interpreted as saying
that trajectories passing near 
infinity have a high cost.
When $\tilde{A}$ has Property \PT,
critical nodes and recurrent nodes coincide (Theorem~\ref{crit-rec}).

Still assuming that $\tilde{A}$ has Property \PT,
we show in Section~\ref{sec-astar-eig}
that each $\rho(A)$-eigenvector satisfying a certain tightness condition
can be represented as an infinite max-plus linear combination
of the critical columns of $\tilde{A}^*$ (Theorem~\ref{generate-astar}).
This extends the previously mentioned theorems of%
~\cite{gondran77,cuning79,cohen83} describing the eigenspace
in the finite dimensional case.
\mytypeout{DISCUSS GONDRAN MINOUX LATEX}

In Section~\ref{sec-cyc}, 
we obtain asymptotic results for the iterates of a matrix
having Property \PT. Theorem~\ref{cycl-th}
extends the cyclicity theorem of~\cite{cohen83}
for finite irreducible matrices. 
Our result is closely related to a ``Turnpike theorem''
that we state as Theorem~\ref{th-turnpike}.

The representation theorem of Section~\ref{sec-astar-eig}
yields only eigenvectors having 
eigenvalue $\rho(A)$, and even for this eigenvalue,
it gives no information about eigenvectors that do
not satisfy the tightness condition.
This is related to the fact that the optimal trajectory associated 
to an eigenvector may head off to infinity. In order to
describe the entire eigenspace, one must develop
a boundary theory analogous to the Martin boundary theory for Markov chains.
We do this in another paper~\cite{AGW-m}. 
In Section~\ref{s-martin} of the present
paper, we recall the principal results without proof, 
and revisit some of the examples of Sections~\ref{sec-rec} 
to~\ref{sec-eigstar} in the light of max-plus Martin boundaries.

\section{Graphs and potential kernels}\label{sec-graph}
Since the supremum of an infinite set may be infinite,
we shall occasionally need to consider the completed
max-plus semiring $\rmaxb$, which is obtained 
by adjoining to $\rmax$ an element $+\infty$,
with the convention that $\zero=-\infty$
remains absorbing for the semiring multiplication.

As in the introduction, we use the symbol $A$ to denote both a matrix 
(or kernel)
$A\in \rmax^{S\times S}$ and the associated max-plus linear
operator $A:\rmax^S\to\rmaxb^S,\; u\mapsto Au$
defined by~\eqref{e-kernel}.
We shall also need to consider infinite valued kernels
and operators, that is kernels $A\in \rmaxb^{S\times S}$
and their corresponding operators $A:\rmaxb^S\to\rmaxb^S,\; u\mapsto Au$.
In addition to being max-plus linear as defined 
in the introduction, such operators are infinitely max-additive:
\begin{align*}
A (\bigoplus_{\ell\in L} u^\ell)&=\bigoplus_{\ell\in L} A u^\ell\quad
\text{for all families }\{u^\ell\}_{\ell\in L}\subset\rmaxb^S\enspace.
\end{align*}

To any matrix $A\in \rmax^{S\times S}$, we associate 
the digraph $G(A)$ with set of nodes
$S$ and an arc $i\to j$ between each pair of nodes $i$ and $j$ 
such that $A_{ij}\neq \zero$. 
The \new{multigraph of strongly connected components} of $G$ is the multigraph
whose nodes are the strongly connected components of $G$,
and which has a number of arcs
between any two such components $C$ and $C'$ equal to the number
of arcs $i\to j$ in $G$, with $i\in C$ and $j\in C'$. 
We say that $G$ is {\em right} (respectively {\em left}) \new{locally finite}
if there are only finitely many
arcs starting from (respectively ending at) each of its nodes.
When $G=G(A)$, we speak about strongly connected components of $A$ instead of strongly connected components of $G(A)$,
etc. We say that $A$ is \new{irreducible} when $A$ has only one 
strongly connected component,
that is, when $G(A)$ is strongly connected. 

Rather than restricting our attention to eigenvectors, we
will also consider the \new{super-eigenvectors}
associated to a value $\lambda$ (or \new{$\lambda$-super-eigenvectors}).
These are vectors
$u\in \rmax^S\setminus\{\zero\}$ such that $Au\leq \lambda u$,
where $\leq$ denotes the pointwise ordering of $\rmax^S$.
The set of all vectors $u\in \rmax^S$ satisfying $Au= \lambda u$
(respectively $Au\leq \lambda u$) is called the \new{$\lambda$-eigenspace}
(respectively the \new{$\lambda$-super-eigenspace}) of $A$. 
We shall say that $u$
has \new{full support} when $u_i\neq \zero$ for all $i\in S$.
Observe that if $A$ is irreducible, then every super-eigenvector
of $A$ has full support.

We define the \new{maximal circuit mean} of $A$ to be
\begin{align*}
 \rho(A) & := \bigoplus_{k\geq 1} (\tr A^k)^{1/k} \in \rmaxb\enspace ,
\end{align*}
where $\tr A:= \bigoplus_{i\in S} A_{ii}$.
Note that $\rho(A)> \zero$ when $A$ is irreducible.
The term ``maximal circuit mean'' can be justified
by giving each path $p$ through $G(A)$
a \new{weight} $|p|_A:=\sum_{(i\to j)\text{ in }p} A_{ij}$
and a \new{length} $|p|$ equal to the number of arcs of $p$. Then, 
\begin{align}
 \rho(A)  := \sup_{c\text{ circuit of } G(A)} \frac{\;|c|_A}{|c|}
\label{e-def-mcm2}
\end{align}
(the division is in the usual algebra). 
The introduction of $\rho(A)$ is motivated
by the following result taken from~\cite[Prop.~3.5]{dudnikov}.
We give a proof for completeness.
\begin{lemma}\label{l1}
If there exists a $\lambda$-super-eigenvector $u$
with full support, then  $\lambda\geq \rho(A)$.
\end{lemma}
\begin{proof}
We have $A^ku\leq \lambda^ku$ and hence 
$(A^k)_{ii} u_i \leq \lambda^k u_i$ for all $i\in S$.
Cancelling $u_i$ and summing over $i\in S$,
we get $\tr A^k \leq \lambda^k$.
Taking the $k$th root and summing over $k\geq 1$,
we get $\rho(A)\leq \lambda$. 
\end{proof}

Given any matrix $A\in \rmax^{S\times S}$, we
define
\begin{align*}
A^+:=A\oplus A^2\oplus A^3 \oplus \cdots\in \rmaxb^{S\times S},
\quad \text{ and }
A^*:=I\oplus A^+ \enspace,
\end{align*}
where $I=A^0$ denotes the max-plus 
identity matrix, that is the matrix with $\unit$
on the diagonal and $\zero$ everywhere else.
We have the following identities
\begin{align*}
A^+=AA^*=A^*A,\quad\text{ and } A^*=A^*A^*\enspace.
\end{align*}
The matrix $A^*$ is analogous to the
\new{potential kernel} in Probabilistic Potential Theory.

It is possible that some entries of $A^*$ take the value $+\infty$.
Indeed, if $\rho(A)>\unit$, then
$A^*_{ii}=+\infty$ for any $i$ belonging
to a circuit of weight greater than $\unit$.
The following proposition, however, gives sufficient conditions
for $A^*$ to be finite.
\begin{prop}[Finiteness condition for $A^*$]
\label{p-fin}
Suppose that $A\in \rmax^{S\times S}$
is such that $\rho(A)\leq \unit$ and the set of paths
between any two distinct nodes of the multigraph of 
strongly connected components
of $A$ is finite. Then, the entries of $A^*$ do not
take the value $+\infty$.
This is the case in particular if $A$ is irreducible.
\end{prop}
\begin{proof}
Since $\rho(A)\le\unit$, we have $A^*_{ii}\leq\unit$ for all $i\in S$.
But $A^*_{ii}=(A^*A^*)_{ii}\geq A^*_{ij}A^*_{ji}$ for all $i,j\in S$, 
and so we get $\unit \geq A^*_{ij} A^*_{ji}$.
When $i$ and $j$ are in the same strongly connected component,
neither $A^*_{ij}$ nor $A^*_{ji}$ equal $\zero$.
Therefore, neither of them can equal $+\infty$.

Consider now two nodes $i$ and $j$ in distinct
strongly connected components, and a path $p$ from $i$ to $j$,
passing successively through the strongly connected components
$C_1,\ldots, C_k$. We can
write $p$ as a concatenation $p=p_1a_1p_2\ldots a_{k-1} p_k$,
where for all $1\leq m\leq k$, $p_m$ is a (possibly empty) path from
some node $i_m\in C_m$ to some node $j_m\in C_m$, 
and for all $1\leq m\leq k-1$, $a_m$ is the arc from $j_m$ to
$i_{m+1}$. Of course, $i_1=i$ and $j_k=j$.
We have
\begin{align*}
|p|_A\leq A^*_{i_1j_1} A_{j_1i_2}
A^*_{i_2j_2} A_{j_2i_3} 
\ldots A_{j_{k-1}i_k}A^*_{i_kj_k}\enspace.
\end{align*}
The assumption on the multigraph
of strongly connected components means that $(i_1,j_1,\ldots,i_k,j_k)$ can only take 
finitely many values. Since
all the $A^*_{i_mj_m}$ are different from $+\infty$, 
we conclude that $A^*_{ij}$, the supremum of $|p|_A$
over all paths from $i$ to $j$, is different from $+\infty$.
\end{proof}
We will denote by $A_{\cdot i}$ and $A_{i \cdot}$, respectively,
the $i$th column and row of a matrix $A$.  
\begin{prop}\label{p-super}
For each $\lambda\in \R$ and $i\in S$,
the column $u=(\lambda^{-1}A)^*_{\cdot i}\in\rmaxb^S$ satisfies
$Au\leq \lambda u$. 
\end{prop}
\begin{proof}
Consider $B:=\lambda^{-1}A$. We have $BB^*=B^+\leq B^*$, 
and taking column $i$, we get $BB^*_{\cdot i}\leq B^*_{\cdot i}$.
In other words $Au\leq \lambda u$.
\end{proof}
\begin{remark}
Combining Propositions~\ref{p-super} and~\ref{p-fin},
we see that when there are only finitely many
paths between any two distinct nodes of the multigraph of strongly connected components
of $A$, any column 
$u=(\lambda^{-1}A)^*_{\cdot i}$,
with $+\infty>\lambda\geq \rho(A)$,
is a super-eigenvector of $A$ associated to the value $\lambda$.
\end{remark}

\section{Recurrent and critical nodes}\label{sec-rec}
{From} now on, we assume that $\rho(A)\in \R$,
and we introduce the normalized matrices
\[
\Al:=\lambda^{-1} A\mbox{ for each $\lambda \in \R$,}\qquad
\text{ and } 
\tilde{A}:=\rho(A)^{-1} A \enspace .
\]
\begin{definition}[Recurrence]\label{defi-rec}
A node $i$ is \new{recurrent}
if $\tilde{A}^+_{ii}=\unit$. We denote by
$N^r(A)$ the set of recurrent nodes.
We call \new{recurrence classes} of $A$ the
equivalence classes of $N^r(A)$ associated to the relation 
$\cR$ defined by $i\cR j$
if $\tilde{A}^+_{ij}\tilde{A}^+_{ji}=\unit$.
\end{definition}
This should be compared with the classical
definition of recurrence for Markov chains,
where a node is recurrent if one
returns to it with probability one.
Here, a node is recurrent
if we can return to it with weight
$\unit$ in $\tilde{A}$.

The notion of recurrence extends the notion
of critical graph occurring in the finite dimensional theory. Recall that
the \new{critical graph} of $A$, denoted by $G^c(A)$, is the subgraph
of $G(A)$ obtained by taking the union of all circuits
of $G(A)$ attaining the maximum
in~\eqref{e-def-mcm2}. The nodes of $G^c(A)$ are
called the \new{critical nodes}, the set of which we denote by $N^c(A)$.
The set of nodes of a strongly connected
component of the critical graph is called
a \new{critical class}.
Obviously, critical nodes are recurrent and
critical classes are contained within recurrence
classes, with equality when $S$ is finite.

When $S$ is infinite, $A$ may have neither critical nodes nor
recurrent nodes.
Note that the definition of recurrence above is
only useful when $S$ is discrete. 
\begin{example}\label{ex-rec-noncrit}
Let $S:=\N$ and take $A_{i,i+1}:=0$ for all $i\in\N$,
$A_{i,0}:=-1/i$ for all $i\in \N\setminus\{0\}$, and 
$A_{ij}:=-\infty$ elsewhere.
The graph of $A$ is:
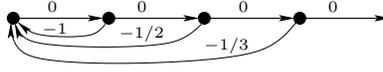
\begin{center}
\input nocritical2
\end{center}
where the nodes are numbered from left to right.
Clearly, $\rho(A)=\unit$. Although $A$ has no critical circuits,
it has a recurrence class, $\N$, since $A^+_{ij}=\unit$ for all
$i,j\in \N$.
\end{example}

\begin{example}\label{ex-rec-noncrit2}
Take the previous example, but now define $A_{i,0}:=-1$
for all $i\in \N\setminus\{0\}$.
The graph of $A$ becomes:
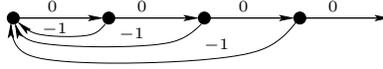
\begin{center}
\input nocritical3
\end{center}
In this example, $\rho(A)=\unit$ and $A$ has no recurrence classes.
\end{example}

\begin{example}\label{ex-rec-noncrit3}
Take again Example~\ref{ex-rec-noncrit}, but replace the diagonal terms
by $A_{ii}=0$ for all $i\in\N$.
The graph of $A$ is now
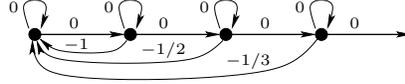
\begin{center}
\input nocritical4
\end{center}
Now $\rho(A)=\unit$ and $A$ has an infinite number 
of critical classes consisting of the singletons $\{i\}$ for
$i\in\N$. Again, $\N$ is a recurrence class.
\end{example}

\section{Tightness and recurrence classes}
We now give a condition which ensures that
the recurrence classes are the same as the critical classes.

We say that a vector $v=(v_j)_{j\in S}\in \rmaxb^S$,
or a map $v:j\mapsto v_j$ from $S$ to $\rmaxb$,
is \new{tight} if for all $j\in S$, $v_j\in\rmax$, 
and if for all $\beta\in \R$, the \new{super-level set}
$J_\beta:=\set{j\in S}{v_j\geq \beta}$ is 
finite. 

This is related to tightness of
idempotent measures
and capacities.
Indeed,
consider the idempotent measure $K_{v}$ with density $v$,
defined by
\[
K_v(J):=\sup_{j\in J} v_j \in \rmaxb \enspace
\qquad\mbox{for all $J\subset S$.}
\]
Then, $v$ is tight in the above sense if and only if
$K_{v}$ is finite and tight in the sense
of~\cite{brienv,12,aqv95,DENSITE,tolya01},
using the discrete topology on $S$.
Tightness conditions also appear in~\cite[Ch.~4]{kim}.

We shall need the following elementary result,
which is a special case of the fact that a function with compact super-level
sets attains its maximum.
\begin{lemma}\label{lem-supattained}
If a map $j\mapsto v_j$ from $S$ to $\rmaxb$ is tight, then the supremum
of $v_j$ over $j\in S$ is attained.\qed
\end{lemma}
We shall also use the following immediate observation.
\begin{obs}\label{obs-tight-leq}
Let $u$ and $v$ be maps from $S$ to $\rmaxb$ and let $\alpha\in\R$.
If $u\leq \alpha v$ and $v$ is tight, then $u$
is tight.\qed
\end{obs}
Let $A\in\rmaxb^{S\times S}$.
We say that 
a vector $u=(u_j)_{j\in S}\in \rmaxb^S$,
or a map $u:j\mapsto u_j$ from $S$ to $\rmaxb$ is
\new{$A$-tight} if, for all $i\in S$,
the map $j\mapsto A_{ij}u_j$, from $S$ to $\rmaxb$, is tight.
\begin{definition}[Property \PT]\label{def-ptyt}
We say that a matrix $A$ \new{has Property \PT} if all the columns
of $A^*$ are $A^*$-tight.
\end{definition}
The following observation shows that this
assumption is symmetric between rows and columns.
\begin{obs}\label{obs-2}
A matrix $A$ has Property \PT\/ if and only if, for all $i,j\in S$, 
the map $k\mapsto A^*_{ik}{A}^*_{kj}$,
from $S$ to $\rmaxb$, is tight.\qed
\end{obs}
\begin{prop}
\label{prop-lequnit}
If a matrix $A$ has Property \PT, then all the entries of ${A}^*$ belong
to $\rmax$, and $\rho(A)\leq \unit$.
\end{prop}
\begin{proof} For all $i,j\in S$, $A^*_{ij}$ is equal to the supremum 
of ${A}^*_{ik}{A}^*_{kj}$ over $k\in S$.
By Observation~\ref{obs-2}, Lemma~\ref{lem-supattained} and the
definition of a tight map, this supremum is attained and must necessarily
belong to $\rmax$.
As already observed before Proposition~\ref{p-fin}, 
this implies that $\rho(A)\leq \unit$. 
\end{proof}
\begin{lemma}\label{atight-irre}
Assume that $A$ is irreducible. Then, the following three statements
are equivalent:
\begin{itemize}
\item[-]
$A$ has Property \PT;
\item[-]
for some $j\in S$, the column ${A}^*_{\cdot j}$ 
is ${A}^*$-tight;
\item[-]
for some $i,j\in S$, 
the map $k\mapsto {A}^*_{ik}{A}^*_{kj}$,
from $S$ to $\rmaxb$, is tight.
\end{itemize}
If these statements are true, then $S$ is denumerable.
\end{lemma}
\begin{proof}
Let $i,j,i',j'\in S$. Since $A$ is irreducible, both
$A^*_{ii'}$ and $A^*_{j'j}$ are strictly greater than $\zero$.
For all $k\in S$, we have
$A^*_{ik}A^*_{kj} \geq A^*_{ii'}A^*_{i'k} A^*_{kj'}A^*_{j'j}$,
and hence, $A^*_{i'k} A^*_{kj'}\leq 
(A^*_{ii'}A^*_{j'j})^{-1} A^*_{ik}A^*_{kj}$.
Using Observation~\ref{obs-tight-leq}, we deduce
that the map $k\mapsto A^*_{i'k} A^*_{kj'}$ is tight whenever the 
map $k\mapsto A^*_{ik} A^*_{kj}$ is tight.
The required equivalences now follow from Observation~\ref{obs-2}.

Fix $i\in S$ and, for each $n\in\N$, let 
$J_{-n}:=\set{k\in S}{{A}^*_{ik}{A}^*_{ki} \geq -n}$.
Since $A$ is irreducible, ${A}^*_{ik}{A}^*_{ki}>\zero$
for all $k\in S$, and so $S=\bigcup_{n\in\N} J_{-n}$.
If $A$ has Property \PT,
then all the sets $J_{-n}$ are finite, in which case,
$S$ will be denumerable.
\end{proof}

\begin{lemma}\label{lem-cycl1}
Let $i,j\in S$ and $\beta>\zero$, and 
let $J:=\set{k\in S}{ A^*_{ik} A^*_{kj}\geq \beta}$. Then,
\begin{align*}
A^{n}_{ij}\leq \beta \oplus (A_{JJ})^n_{ij} \quad\text{for all } 
n\geq 0 \enspace .
\end{align*}
\end{lemma}
\begin{proof} 
Let $p:=(i_0=i,\ldots, i_n=j)$ be a path from $i$ to $j$ of length $n$.
 If all the nodes of $p$ are contained in
$J$, then $|p|_A\leq (A_{JJ})^n_{ij}$. If not, then one of the nodes $i_q$,
with $0\leq q\leq n$, is in $S\setminus J$, in which case
$|p|_A\leq A^q_{i,i_q} A^{n-q}_{i_q,j} \leq  A^*_{i,i_q} A^*_{i_q,j}
< \beta $.
\end{proof}

\begin{prop}\label{lem-eq16}
Assume that $A$ has Property \PT.
Let $i,j\in S$ and $n\in\N$ be such that $A^n_{ij}\neq\zero$.
Then, there exists a path $p:=(i_0=i,i_1,\ldots, i_n=j)$ 
of length $n$ from $i$ to $j$ such that 
\begin{align*}
A^n_{ij}=|p|_A:=A_{i_0 i_1}\cdots A_{i_{n-1} i_n}\enspace . \end{align*}
Moreover, if $A^+_{ij}\neq \zero$, then  there exists an elementary path
$p$ of nonzero length from $i$ to $j$ such that $A^+_{ij}=|p|_A$.
\end{prop}
\begin{proof} Fix $i,j\in S$ and $n\in\N$.
Choose $\beta< A^n_{ij}$ and denote 
$J:=\set{k\in S}{ A^*_{ik} A^*_{kj}\geq \beta}$. 
Then, $i,j\in J$ and, by Property \PT, $J$ is a finite set.
By Lemma~\ref{lem-cycl1},
$A^{n}_{ij}\leq \beta \oplus (A_{JJ})^n_{ij}$.
But, $A^{n}_{ij}>\beta$ and
$A^{n}_{ij}\geq (A_{JJ})^n_{ij}$,  so
$A^{n}_{ij}= (A_{JJ})^n_{ij}$.
Now, there are only a finite number of paths
of length $n$ from $i$ to $j$ that are contained entirely within $J$.
Hence, the supremum of $|p|_A$ over the set of these
paths is attained by some path
$p$, and for this path $(A_{JJ})^n_{ij}=|p|_A$.

Choose now $\beta'< A^+_{ij}$ and let
$J:=\set{k\in S}{ A^*_{ik} A^*_{kj}\geq \beta'}$. 
Again, $i,j\in J$ and $J$ is finite.
Applying Lemma~\ref{lem-cycl1} and taking the max-plus sum over
$n\geq 1$, we obtain that
$A^{+}_{ij}\leq \beta' \oplus (A_{JJ})^+_{ij} $.
This implies that $A^{+}_{ij}=(A_{JJ})^+_{ij}$.
Now, any path $p$ from $i$ to $j$ can be decomposed into a
disjoint union of an elementary path $p_0$ from $i$ to $j$ and elementary
circuits. Since, by Proposition~\ref{prop-lequnit}, $\rho(A)\leq \unit$,
we have that $|p|_A\leq |p_0|_A$. Hence $(A_{JJ})^+_{ij}$ is
the supremum of $|p|_A$ over all elementary paths $p$ of length $|p|\geq 1$
from $i$ to $j$ in $J$. But since $J$ is finite, there are only finitely
many such paths and so the supremum is attained.
\end{proof}
\begin{theorem} \label{crit-rec}
Assume that $\tilde{A}$ has Property \PT.
Then $N^c(A)=N^r(A)$ and the critical classes coincide with the recurrence
classes. Furthermore, the critical graph $G^c(A)$ coincides with the graph
$G^r(A)$ having set of nodes $N^r(A)$ and an arc $(i,j)$ whenever
$\tilde{A}_{ij}\tilde{A}^+_{ji}=\unit$.
\end{theorem}
\begin{proof} 
We may assume that $\rho(A)=\unit$,
so that $\tilde{A}=A$.
Let us first prove that $N^r(A)=N^c(A)$.
Recall that each critical node is recurrent, that is 
$N^c(A)\subset N^r(A)$. To prove the reverse inclusion,
suppose $i\in S$ is recurrent, in which case $A^+_{ii}=\unit$.
By the second part of Proposition~\ref{lem-eq16},
there exists an elementary circuit containing
$i$ with weight $\unit$. But $\rho(A)= \unit$, 
and so $i$ is critical.

Let $i$ and $j$ be nodes in the same recurrence class, so that
$A^+_{ij}A^+_{ji}=\unit$.
Applying Proposition~\ref{lem-eq16} again, we see that
there exists an elementary path $p$ from $i$ to $j$ 
such that $A^+_{ij}=|p|_A$ and an elementary path $p'$
from $j$ to $i$ such that $A^+_{ji}=|p'|_A$.
Concatenating these two paths, we obtain a critical circuit 
passing through both $i$ and $j$. Hence, $i$ and $j$ are in the same critical 
class.  This shows that each recurrence class is contained within a critical
class. Since the converse inclusion is trivial,
the recurrence and critical classes must coincide.

Now, let $(i,j)$ be an arc of $G^r(A)$, so that
$A_{ij} A^+_{ji}=\unit$.
Again, by Proposition~\ref{lem-eq16} 
there exists an elementary path $p$ from $j$ to $i$ 
such that $A^+_{ji}=|p|_A$. Concatenating this path with the arc $(i,j)$,
we obtain a critical  circuit containing $(i,j)$.
Hence this arc is an arc of $G^c(A)$.
The converse, that each
arc of the critical graph $G^c(A)$ is an
arc of $G^r(A)$, is trivial.
\end{proof}
None of the matrices in the three examples~\ref{ex-rec-noncrit}--\ref{ex-rec-noncrit3} have Property \PT.
Nevertheless, the conclusion of Theorem~\ref{crit-rec} is true for
Example~\ref{ex-rec-noncrit2}.
\begin{example}\label{ex-assum-tight1}
Let $S:=\N$ and take $A_{i,i+1}:=0$ for all $i\in\N$,
$A_{00}:=0$, $A_{i,i-1}:=-1$ for all $i\in \N\setminus\{0\}$, and 
$A_{ij}:=-\infty$ elsewhere.
The graph of $A$ is
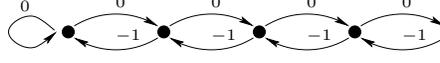
\begin{center}
\input tight1
\end{center}
Obviously, $A$ is irreducible and $\rho(A)=\unit$.
Calculating, we find that $A^+_{ij}$ takes the value $j-i$ when $i>j$,
the value $-1$ when $i=j\neq 0$, and the value $0$ otherwise.
So, $A$ has Property \PT. We observe that 
$A$ has a single critical class $\{0\}$,
which of course is also the only recurrence class.
\end{example}
\begin{example}\label{ex-assum-tight2}
Let $S:=\N$ and take $A_{i,i+1}:=0$ for all $i\in\N$,
$A_{i,i-1}:=-1/i$ for all $i\in \N\setminus\{0\}$, and 
$A_{ij}:=-\infty$ elsewhere.
The graph of $A$ is
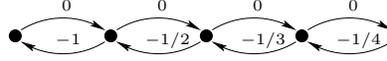
\begin{center}
\input tight2
\end{center}
Again, $A$ is irreducible and
$\rho(A)=\unit$. This time $A^+_{ij}=-\sum_{k={j+1}}^i 1/k$ when $i>j$,
$A^+_{ij}=-1/(i+1)$ when $i=j$, and $A^+_{ij}=0$ otherwise.
So this $A$ also satisfies Property \PT.
However, now there are no critical classes and no recurrence classes.
\end{example}
\section{Eigenvectors associated to recurrent nodes}\label{sec-eigstar}
The following proposition extends a well known result
in the finite dimensional case.
\begin{prop}[Recurrent columns of $\tilde{A}^*$ are eigenvectors]\label{l2}
Let  $\lambda\in\R$ be such that $\lambda\geq \rho(A)$.
Then, the vector $u:=(\Al )^*_{\cdot i} $ satisfies
$\lambda u=Au$ if and only if $i$ is a recurrent node
and $\lambda=\rho(A)$.
In particular, if $\rho(A)\in\R$ and there are only finitely many paths
between any two distinct nodes of the multigraph of strongly connected components of $A$,
then $\tilde{A}^*_{\cdot i}$ is a
$\rho(A)$-eigenvector of $A$
if and only if $i$ is recurrent.
\end{prop}
\begin{proof}
We have $A(\Al )^*=\lambda (\Al )^+$.
Observe that
$(\Al )^*$  and $(\Al )^+$ coincide, 
except perhaps on the diagonal,
where $(\Al )^*_{ii}=\unit \oplus (\Al )^+_{ii}$.
But if $i$ is recurrent and $\lambda=\rho(A)$,
then $(\Al )^+_{ii}= \tilde{A}^+_{ii}=\unit$,
and so the $i$th columns of $(\Al )^*=\tilde{A}^*$
and $(\Al )^+=\tilde{A}^+$ coincide.
In this case, 
\[ 
A({A}_\lambda)^*_{\cdot i}=\rho(A) \tilde{A}^+_{\cdot i}
   =\rho(A) \tilde{A}^*_{\cdot i}=\rho(A) ({A}_\lambda)^*_{\cdot i}
\enspace .\]
Conversely, if $A(\Al )^*_{\cdot i}=\lambda  (\Al )^*_{\cdot i}$,
then \[(\Al )^+_{ii}=(\Al  (\Al )^*_{\cdot i})_{i}
= (\Al )^*_{ii}=\unit\enspace .\]
This implies that $\lambda=\rho(A)$ and $i$ is recurrent.

We now assume that there are only finitely many paths
between any two distinct nodes of the multigraph of strongly connected components of $A$.
The only additional thing to prove in this case is that 
$u:={\tilde{A}}^*_{\cdot i}\in\rmax^S\setminus\{\zero\}$.
But this follows by applying Proposition~\ref{p-fin} to $\tilde A$.
\end{proof}
The following observation shows that eigenvectors
corresponding to nodes in the same recurrence class are proportional.
\begin{prop}\label{l22}
If $i$ and $j$ belong to the same recurrence class, 
then 
\[\tilde{A}^*_{\cdot i}=\tilde{A}^*_{\cdot j}\tilde{A}^*_{ji}\enspace. \]
\end{prop}
\begin{proof}
For $i$ and $j$ belonging to the same recurrence class,
we have $\tilde{A}^+_{ij}\tilde{A}^+_{ji}=\unit$,
and thus $\tilde{A}^*_{ij}\tilde{A}^*_{ji}=\unit$.
Using this and the fact that $(\tilde{A}^*)^2=\tilde{A}^*$,
we see that
\[\tilde A^*_{ki}\geq \tilde A^*_{kj}\tilde A^*_{ji}
\geq \tilde A^*_{ki}\tilde A^*_{ij}\tilde A^*_{ji}
= \tilde A^*_{ki}\]
for all $k\in S$.
Therefore
$ \tilde A^*_{ki}= \tilde A^*_{k j}\tilde A^*_{ji}$ for all 
$k\in S$.
\end{proof}
\begin{example}\label{ex-birth}
Let $S:=\N$, $p,q\in \rmax$ and take
$A_{i,i+1}:=p$ for $i\in\N$, $A_{i,i-1}:=q$
for $i\in\N\setminus\{0\}$, 
and $A_{ij}:=\zero$ elsewhere. 
Consider first the case when neither $p$ nor $q$ are equal to $\zero$.
Then, the graph of $A$ is:
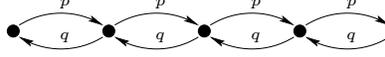
\begin{center}
\input birth2
\end{center}
We have $\rho(A)=\sqrt{pq}$, that is $\rho(A)=(p+q)/2$ in the
usual algebra notation. 
Since $A$ is irreducible,
any eigenvector of $A$ has full support.
Lemma~\ref{l1} tells us that any eigenvalue $\lambda$ must be greater
than or equal to $\rho(A)$. Since $\N$ is a critical class,
therefore necessarily a recurrence class, 
Proposition~\ref{l2} shows that
every column of $\tilde{A}^*$ is an eigenvector
of $A$ with eigenvalue $\rho(A)$ and Proposition~\ref{l22}
shows that they are all proportional. 
We readily check that for any $\lambda \geq \rho(A)$, 
the vector $u\in\R^S$ defined by $u_k :=(\lambda/p)^k$ (that is,
$u_k:=k(\lambda -p)$ in the usual algebra notation)
is an eigenvector with eigenvalue $\lambda$.
So the \new{spectrum} of $A$,
that is, the set of eigenvalues of $A$,
is $[\rho(A),+\infty)$.

Now consider the case when $pq=\zero$. This time $\rho(A)=\zero$ and $A$
is not irreducible.
When $p\neq \zero$ and $q=\zero$, the vectors given by
$u_k :=(\lambda/p)^k; k\in\N$ are still
eigenvectors, and so the
spectrum of $A$ is $[\zero,+\infty)$. 
When $p=\zero$ and $q\neq \zero$, the spectrum of $A$ is $\{\zero\}$.

Since transposing $A$ corresponds to exchanging $p$ and $q$,
we deduce that $A$ and its transpose
have the same spectrum when $p\neq\zero$ and $q\neq\zero$. 
\end{example}
\begin{example}\label{ex-triang}
We now give an example to show that the spectrum of an irreducible
matrix may differ from that of its transpose.

Consider $S=\N$ and the matrix
\[
A=\begin{pmatrix}
\alpha_0   & \beta          &  \beta  & \ldots \\
\beta         & \alpha_1    & \beta     & \dots\\
\zero      & \beta          & \alpha_2       & \ddots\\
\vdots      &  \ddots  & \ddots &\ddots 
\end{pmatrix}\enspace.
\]
We assume that $\zero<\beta<\unit$, that
$\zero<\alpha_i<\unit$ for all $i\in\N$, and that 
$\lim_{i\to\infty} \alpha_i =\unit$. Under these assumptions,
$\rho(A)=\unit$, $A$ is irreducible,
and the entries of $\tilde{A}^*=A^*$ are finite. 
Observe that the maximal circuit mean $\rho(A)$ is obtained by taking
a loop at node $i$ and letting $i$ tend to $+\infty$.
There are no critical nodes nor recurrent nodes.

We claim that the spectrum of $A$ is empty.
Indeed, suppose $u$ is an eigenvector with eigenvalue $\lambda$.
Then, $u$ has full support since $A$ is irreducible,
and Lemma~\ref{l1} gives that $\lambda\geq \rho(A)=\unit$.
Setting $\mu=\lambda^{-1} \beta$ and $\gamma_k=\lambda^{-1}\alpha_k$,
we rewrite the spectral equation $Au=\lambda u$
as
\(u_i=\gamma_i u_i \oplus  \mu (\bigoplus_{j\in \N,j\neq i,j\geq i-1} u_j) 
\).
Since $\gamma_i<\unit$, we see that
\begin{equation}\label{uiM}
u_i =  \mu(\bigoplus_{j\in \N, j\neq i,j\geq i-1} u_j)\enspace.\end{equation}
Taking the sum of these equalities over $i\in \N$,
and setting $M:=\bigoplus_{j\in \N}u_j$, we get that $M=\mu M$.
Since $\mu<\unit$, this implies that either $M=\zero$ or $M=+\infty$.
The case where $M=\zero$ can be eliminated, since
$u$ has at least one entry different from $\zero$.
The case where $M=+\infty$ can also be eliminated
since Equation~\eqref{uiM} yields $u_0\geq \mu M$.
It follows that the spectrum of $A$ is empty.

The transpose of $A$, on the other hand, has a non-empty spectrum:
any $\lambda \geq \unit$ is an eigenvalue and the 
corresponding eigenvector is $(\mu^0,\mu^{-1},\mu^{-2},\ldots)$.
\end{example}

We conclude this section with two more properties of super-eigenvectors.
The finite dimensional versions of these were
instrumental in~\cite{ABG96,gg}. 
\begin{prop}
If $\rho(A)\in\R$ and $v,w\in \rmax^S$ are super-eigenvectors
of $A$ associated to the value $\rho(A)$, 
then the restrictions of $v$ and $w$ to
any recurrence class of $A$ are proportional.
\end{prop}
\begin{proof}
Let $v$ be a super-eigenvector of $A$ associated to $\rho(A)$.
Then $\tilde{A}v\leq v$, and so $\tilde{A}^* v\leq v$.
Fix a recurrence class $C$ of $A$.
For all $i,j\in C$, we have  
$v_i\geq \tilde{A}^*_{ij} v_j\geq \tilde{A}^*_{ij} \tilde{A}^*_{ji}
v_i=v_i$, hence $v_i=\tilde{A}^*_{ij} v_j$. 
Therefore, the restriction of $v$ to $C$ is
proportional to the restriction of $\tilde{A}^*_{.j}$ to $C$.
Hence the restriction of any two super-eigenvectors $v$ and $w$ to $C$
are proportional.
\end{proof}

The probabilistic analogue of the next result is related to the
minimum principle, as explained in~\cite[Lemma~2.9 and Th.~3.4]{conv},
in a slightly different context.
\begin{lemma}
If $\rho(A)\in\rmax$, $Au\leq \rho(A)u$ and $u\in \rmax^S$, 
then $(Au)_i=\rho(A)u_i$ for each recurrent node $i$. 
\end{lemma}
\begin{proof}
Assume that $(Au)_i<\rho(A)u_i$ for some $i\in S$.
Then, $\rho(A)\in \R$, $u_i\neq \zero$, and
$(Au)_i\leq \alpha\rho(A)u_i$
for some $\alpha<\unit$. 
For all circuits $(i_1,\ldots,i_k)$ in $G(A)$
starting at $i_1=i$, we have
\[
\tilde A_{i_1i_2}u_{i_2}\leq \alpha u_{i_1},\;\;
\tilde A_{i_2i_3}u_{i_3}\leq  u_{i_2},\;\;
\ldots,\;\;
\tilde A_{i_ki_1}u_{i_1}\leq  u_{i_k}
\enspace.
\]
It follows
that $u_{i_2},\ldots,u_{i_k}$ are also different from $\zero$.
Combining the inequalities, we get that
$\tilde A_{i_1i_2}\cdots \tilde A_{i_ki_1}\leq \alpha$.
Since this holds for all sequences $(i_1,\ldots,i_k)$
such that $i_1=i$, we get that $\tilde A^+_{ii}\leq \alpha<\unit$,
and so $i$ is not recurrent.
\end{proof}

\section{Representation of tight eigenvectors}\label{sec-astar-eig}
We show here that, when looking for $\tilde{A}^*$-tight eigenvectors,
the situation is similar to the case when $S$ is finite.
In particular, the eigenspace can be described in terms of the critical graph.
\begin{theorem}\label{th-unique}
Let $\lambda\in \R$ be such that $\lambda\geq \rho(A)$.
If $u\in\rmax^S\setminus\{\zero\}$ is a $(\Al )^*$-tight 
$\lambda$-eigenvector, then
$\lambda=\rho(A)$, the set $N^c(A)$ of critical nodes of $A$ is non-empty, and
\begin{align}\label{e-rep}
u= \bigoplus_{j\in N^c(A)} \tilde{A}^*_{\cdot j} u_j \enspace .
\end{align}
\end{theorem}
\begin{proof}
We shall assume, without loss of generality, that $\lambda=\unit$,
in which case $\Al =A$.
Let $u$ be a $\unit$-eigenvector.
Then, $u=Au$ and so $u=A^*u$.
It follows that
\begin{align*}
u_i\geq \bigoplus_{j\in N^c(A)} A^*_{ij} u_j\quad \mbox{for all $i\in S$}
\enspace.
\end{align*}
We now show the opposite inequality when $u$ is $A^*$-tight.
If $u_i=\zero$, this is trivial and so we assume the opposite.
Since  $u$ is $A^*$-tight, it is also $A$-tight.
Therefore, by Lemma~\ref{lem-supattained}, for each $j\in S$, there is some
$k\in S$ such that $u_j=A_{jk}u_k$. 
This allows us to construct a sequence
$\{i_k\}_{k\geq 1}\subset S$ such
that $i_1=i$ and $u_{i_k}=A_{i_ki_{k+1}}u_{i_{k+1}}$ for all $k\geq 1$.
For this sequence,
\begin{align}
u_{i_\ell} = A_{i_\ell i_{\ell+1}}\cdots A_{i_{k-1} i_{k}}u_{i_{k}}
\quad\mbox{for all $1\leq \ell\leq k$} \enspace .
\label{e-path}
\end{align}
Suppose that all the $i_k$ are distinct.
Then, the $A^*$-tightness of $u$ implies that, for each $\alpha\in \R$,
the set $J:=\set{j\in S}{A^*_{ij} u_j \geq \alpha}$ is finite.
Now, some $i_{k}$ must be in $S\setminus J$.
So, we may combine this property with~\eqref{e-path} to get
\begin{align*}
u_i \leq A^*_{ii_{k}}u_{i_{k}}\leq \alpha \enspace.
\end{align*}
This holds for all $\alpha\in \R$, and so
$u_i=\zero$, which contradicts a previous assumption.
This shows that the $i_k$ can not be all distinct,
in other words that $i_{k+c}=i_k$ for some $k,c\geq 1$.
{From}~\eqref{e-path} with $\ell=1$, we deduce that $u_{i_k}\neq \zero$.
Using this and Equation~\eqref{e-path} again,
we get that
$A_{i_ki_{k+1}}\cdots A_{i_{k+c-1}i_k}=\unit$.
This shows that $\rho(A)\geq \unit$, and since we have assumed that
$\rho(A)\leq \lambda=\unit$, we obtain that $\rho(A)= \unit$.
Moreover, $i_k\in N^c(A)$.
Using~\eqref{e-path} again with $\ell=1$, we get
\begin{align*}
u_i &\leq A^*_{ii_k} u_{i_k} \leq
\bigoplus_{j\in N^c(A)}A^*_{ij} u_{j}.
\typeout{BEHAVIOR MAY DEPEND ON THE VERSION OF AMSTEX, POSITION OF qed TO BE CHECKED CAREFULLY}
\qedhere
\end{align*}
\end{proof}

\begin{corollary}
If there exists an $\tilde{A}^*$-tight 
$\rho(A)$-eigenvector, then there exists an  $\tilde{A}^*$-tight
critical column, that is, a column $\tilde{A}^*_{\cdot j}$, with $j\in N^c(A)$.
\end{corollary}
\begin{proof}
Let $u$ be an $\tilde{A}^*$-tight eigenvector.
By Theorem~\ref{th-unique}, we can find $j\in N^c(A)$ such that
$u_j\neq \zero$. Then $\tilde{A}^*_{\cdot j}\le (u_j)^{-1} u$
and so, by Observation~\ref{obs-tight-leq},
 $\tilde{A}^*_{\cdot j}$ is necessarily $\tilde{A}^*$-tight.
\end{proof}
\begin{corollary}\sloppy
If $A$ is irreducible and there exists an $\tilde{A}^*$-tight 
$\rho(A)$-eigenvector, then $\tilde A$ has Property \PT.
\end{corollary}
\begin{proof}
We apply Lemma~\ref{atight-irre}.
\end{proof}
The set $\sT$ of $\tilde{A}^*$-tight vectors $u$ such that
$Au=\rho(A)u$ is a subsemimodule of $\rmax^S$, meaning 
that $\sT$ is stable under finite max-plus linear combinations.
It has the following additional property:
the least upper bound in $\rmax^S$ of 
any family of elements of $\sT$ that
is bounded above by an element of $\sT$ belongs to $\sT$.
Subsets of $\rmax^S$ 
with the two properties above are said to be \new{boundedly complete
subsemimodules} of $\rmax^S$.
See~\cite{litvinov00,litvinov02,cgq02} for general definitions
concerning semimodules over idempotent semirings.

Given a semimodule $V$ over $\rmax$, one is naturally interested in
representing the elements of $V$ in terms of generators. 
We say that $\xi\in V\setminus\{\zero\}$ is an \new{extremal generator} 
of $V$ if $\xi=u\oplus v$ with $u,v\in V$
implies that either $\xi=u$ or $\xi=v$.
This concept has, of course, an analogue in the usual algebra, where
extremal generators are defined for cones.
Max-plus extremal generators are also called \new{join irreducible} elements
in the lattice literature.
When $V$ is a boundedly complete subsemimodule of $\rmax^S$,
we shall say that 
$V$ is \new{spanned} by a family $\{\xi^k\}_{k\in K}\subset V$ if
any element $v$ of $V$ can be expressed as a (possibly infinite)
max-plus linear combination
of the $\xi^k$, in other words, has a representation
$v=\bigoplus_{k\in K} \nu_k \xi^k$,
where $\nu_k\in\rmax$ for each $k\in K$.
We say that a spanning family of $V$
is \new{minimal} if it does not include any smaller
spanning family of $V$. 
\begin{prop} \label{prop-min-isextremal}
Let $V$ be a boundedly complete subsemimodule of $\rmax^S$, and let
$\{\xi^k\}_{k\in K}$ be a minimal spanning family of $V$.
Then, $\xi^k$ is an extremal generator of $V$
for all $k\in K$.
\end{prop}
\begin{proof}
Assume that $\{\xi^k\}_{k\in K}$ is a spanning family of $V$, but
not necessarily minimal, and suppose it contains
an element $\xi^\ell$ that is not an extremal generator of $V$.
We may write  $\xi^\ell=u\oplus v$ with both $u$ and $v$ in $V$ different
from $\xi^\ell$.
Since $\{\xi^k\}_{k\in K}$ is a spanning family of $V$, 
there exist $\mu=(\mu_k)_{k\in K}\in \rmax^K$
and $\nu=(\nu_k)_{k\in K}\in \rmax^K$ such that
$u=\bigoplus_{k\in K} \mu_k \xi^k$ and $v=\bigoplus_{k\in K} \nu_k \xi^k$.
This implies for instance that $\xi^\ell\geq u \geq \mu_\ell \xi^\ell$.
Hence $\mu_\ell \xi^\ell_i \leq \xi^\ell_i$ for all $i\in S$, and 
since $\xi^\ell$ is different from $\zero$, we deduce that 
$\mu_\ell\leq \unit$. If $\mu_\ell$ was equal to $\unit$, we would have
$\xi^\ell=u$, a contradiction. Therefore, 
$\mu_\ell<\unit$, and similarly $\nu_\ell<\unit$, 
hence $\mu_\ell\oplus\nu_\ell <\unit$.
Since $\xi^\ell=u\oplus v=\bigoplus_{k\in K} (\mu_k\oplus \nu_k) \xi^k$,
and for all $i\in K$, the maximum of $(\mu_k\oplus \nu_k)  \xi^k_i$ 
over $k\in K$ cannot be attained at $k=l$ unless $\xi^\ell_i=\zero$,
we have the representation  
\begin{equation*}
\xi^\ell= \bigoplus_{k\in K\setminus\{\ell\}}
 (\mu_k\oplus \nu_k) \xi^k\enspace.
\end{equation*}
Since $\xi^\ell$ can be expressed as a max-plus linear combination
of $(\xi^k)_{k\in K\setminus\{\ell\} }$, this family must also be
a spanning family of $V$. Thus the  spanning family $\{\xi^k\}_{k\in K}$
is non-minimal.
\end{proof}
The following corollary of Theorem~\ref{th-unique} extends
a basic result of finite dimensional max-plus spectral 
theory~\cite{gondran77,cuning79,cohen83} (see also~\cite[Th.~3.100]{bcoq}).
\begin{theorem}\label{generate-astar}
Assume that $\tilde A$ has Property \PT.
Let $\sT$ be the semimodule of $\tilde{A}^*$-tight vectors
$u\in\rmax^S$ satisfying $Au=\rho(A) u$.
Then, the extremal generators of $\sT$ are precisely
the vectors of the form $\alpha \tilde{A}^*_{\cdot j}$,
where $\alpha\in \R$ and $j$ is a critical node.
Moreover, we may obtain a minimal spanning family of $\sT$ by taking
exactly one column $\tilde{A}^*_{\cdot j}$ from each
critical class of $A$.
\end{theorem}
\begin{proof}
Again we can assume, without loss of generality,
that $\rho(A)=\unit$, and so $\tilde{A}=A$.
By Theorem~\ref{crit-rec}, $N^c(A)=N^r(A)$.
If there are no critical nodes, then by 
Theorem~\ref{th-unique}, $\sT=\{\zero\}$, which is spanned by the
empty set. So we shall assume that $N^c(A)\neq\emptyset$.
Denote by $J$ some subset of $N^c(A)$ obtained by taking
exactly one element from each critical class of $A$.
Any $A^*_{\cdot j}$ with $j\in J$ is a $\rho(A)$-eigenvector
by Proposition~\ref{l2},
and is $A^*$-tight since $A$ has Property \PT.
Thus $A^*_{\cdot j}\in \sT$ for all $j\in J$.

By Proposition~\ref{l22}, all the $A^*_{\cdot j}$ with $j$ in the same
critical class are proportional.
So, the representation of $u\in \sT$ given by
Theorem~\ref{th-unique} can be reduced to
$u= \bigoplus_{j\in J} A^*_{\cdot j} \nu_j$,
where $\nu_j\in \rmaxb$ for all $j\in J$. 
We may rule out that $\nu_j=+\infty$ for any $j\in J$
since this would lead to $u$ having a component $+\infty$ (by
$u_j\geq A^*_{jj} \nu_j=\nu_j$) which is impossible since 
$\sT\subset \rmax^S$.
This shows that $\sT$ is spanned by the family $(A^*_{\cdot j})_{j\in J}$.

We now show that this family is minimal.
Suppose, for some $i\in J$,
$A^*_{\cdot i}= \bigoplus_{j\in J\setminus\{i\}} \nu_j A^*_{\cdot j}$.
Then, 
$A^*_{ji}\geq \nu_j A^*_{jj}=\nu_j$ for all $j\in J\setminus\{i\}$.
So the vector $\nu$ is $A^*$-tight, if we define $\nu_k:=\zero$
for all $k\not\in J\setminus\{i\}$.
Hence, by Lemma~\ref{lem-supattained}, for all $k\in S$, the supremum in
 $A^*_{ki}= \bigoplus_{j\in J\setminus\{i\}} \nu_j A^*_{k j}$
is attained by some $j\in J\setminus\{i\}$, that is
$A^*_{ki}= \nu_j A^*_{k j}$.
Applying this to $k=i$, we get
$\unit=A^*_{ii}=\nu_j A^*_{ij}$.
Since $\nu_j \leq A^*_{ji}$, we obtain 
$\unit\leq  A^*_{ji} A^*_{ij}\leq \unit$.
Hence, $i$ and $j$ are in the same recurrence class, 
and so by Theorem~\ref{crit-rec}, they are in the same critical class,
a contradiction. The second statement of the theorem follows.

We deduce from the second statement of the theorem and from
Proposition~\ref{prop-min-isextremal} that any vector
proportional to some $A^*_{\cdot j}$, where
$j$ is critical, is an extremal generator of $\sT$.
It remains to check that all the extremal generators of $\sT$
are obtained in this way. So, suppose there existed
an extremal generator $u$ of $\sT$ which was not proportional to any of these
$A^*_{\cdot j}$. 
For any $k\in N^c(A)$, we may write the representation given 
in Theorem~\ref{th-unique} in the form $u=u_{k}  A^*_{\cdot k}
\oplus v$, where $v=\bigoplus_{j\in N^c(A)\setminus\{k\}}
u_j A^*_{\cdot j}$. Since $u$ is extremal
and not proportional to $A^*_{\cdot k}$, 
we must have $u=v$. More generally, we may write
\begin{align}
\label{e-outK}
u=\bigoplus_{j\in N^c(A)\setminus K}u_j A^*_{\cdot j}\enspace,
\end{align}
for any finite subset $K$ of $N^c(A)$. 
Consider now any $i\in S$ and $\beta\in \R$, and
let $K:=\set{j\in N^c(A)}{A^*_{ij} u_j\geq \beta}$.
Since $u$ is $A^*$-tight, the set $K$ is
finite, and we deduce from~\eqref{e-outK}
that $u_i\leq \beta$. Since this holds for all $i\in S$ and
$\beta\in\R$, it follows that $u=\zero$,
which contradicts the assumption that it is
an extremal generator of $\sT$.
\end{proof}
We shall see in Section~\ref{s-martin}
that if $\tilde{A}$ does not have Property \PT, then the set
of $\rho(A)$-eigenvectors of $A$
may or may not be spanned by the critical columns of $\tilde{A}^*$.
We shall also see that, even when $\tilde A$ has Property \PT,
there may exist $\rho(A)$-eigenvectors that 
cannot be represented in the form~\eqref{e-rep}.
These, of course, are not $\tilde{A}^*$-tight.

\section{Cyclicity theorem}\label{sec-cyc}
In this section, we investigate the powers
of matrices having Property \PT.
Our first result is a max-plus
analogue of the fact that the powers of a transient
Markov matrix converge to zero.
\begin{theorem}\label{cor-cycl1}
Assume $A$ has Property \PT, and that either
$\rho(A)<\unit$ or $N^c(A)=\emptyset$.
Then, for all $i,j\in S$, 
\[
\lim_{n\to\infty} {{A}}^n_{ij}=\zero \enspace .
\]
\end{theorem}
\begin{proof}
Let us fix $i,j\in S$ and $\beta>\zero$, and 
let $J:=\set{k\in S}{ A^*_{ik} A^*_{kj}\geq \beta}$. 
Since $A$ has Property \PT, $J$ is a finite set.
By Proposition~\ref{prop-lequnit}, $\rho(A)\leq\unit$.
Since either $\rho(A)<\unit$ or $N^c(A)$ is empty, all 
circuits have a weight strictly less than
$\unit$. Since $J$ is finite, we conclude that $\rho(A_{JJ})<\unit$. 
This implies that $(A_{JJ})_{ij}^n$ converges 
to $\zero$ as $n\to\infty$. Indeed,
for any matrix $B\in \rmax^{d\times d}$
either $\rho(B)=\zero$,
and then $B^n=\zero$ for all $n\geq d$
(because any path of length at least $d$ in the graph
of $B$ must contain a circuit), or $\zero<\rho(B)$,
and then, $\tilde{B}^*$ exists and is finite, and
$B^n\leq \rho(B)^n\tilde{B}^*$ for all $n\geq 0$. 
Taking $B=A_{JJ}$, we use Lemma~\ref{lem-cycl1} to deduce that 
$A^{n}_{ij}\leq \beta $ for $n$ large enough,
and the result follows since $\beta$ is arbitrary.
\end{proof}
\begin{remark}
The convergence of $A^n_{ij}$ to $\zero$ may be arbitrary slow.
For instance the matrix
$A$ of Example~\ref{ex-assum-tight2}
satisfies the assumptions of Theorem~\ref{cor-cycl1},
and $A_{00}^{2n}=-\sum_{k=1}^n \frac{1}{k}\sim-\log n$,
which goes to $\zero$ in a sublinear way.
\end{remark}
\begin{definition}[Cyclicity]
The \new{cyclicity} $\cyc{G}$ of a strongly connected (possibly infinite)
graph $G$ is defined to be the $\gcd$ of the lengths of its circuits.
The cyclicity of a non strongly connected graph is defined to be the
$\lcm$ of the cyclicities of its strongly connected components.
\end{definition}
For a matrix $A\in \rmax^{S\times S}$,
we denote by $\cyc{A}$ the cyclicity of its graph $G(A)$ and
by $\cycc{A}$ the cyclicity of its critical graph $G^c(A)$.

Note that $\cyc{G}$ may be infinite when $G$ has an infinite number of
strongly connected components. If there are no strongly connected components,
then $\cyc{G}$ is equal to $1$
(by convention, the $\lcm$ of an empty set is $1$).
We write $a\=mod b[c]$ when $a-b$ is a multiple of $c$.
\begin{theorem}[Cyclicity theorem]\label{cycl-th}
Assume that $A$ is irreducible, that $\tilde A$ has Property \PT,
and that $N^c(A)\neq\emptyset$. Then, for all $i,j\in S$, there exist
$\sigma_{ij}\in\N\setminus\{0\}$ and 
 $n_{ij}\in\N$ such that,
\begin{equation}\label{cycl1p3}
A^{n+\sigma_{ij} }_{ij}=\rho(A)^{\sigma_{ij} } A^n_{ij}
\quad\text{ for }n\geq n_{ij}\enspace,
\end{equation}
and we have the following explicit formula
\begin{equation}\label{cycl1}
\tilde{A}^{n}_{ij}=\bigoplus_{k\in N^c(A)} 
(\tilde{A}^q (\tilde{A}^{\sigma_{ij} })^*)_{ik} 
(\tilde A^{q'}(\tilde{A}^{\sigma_{ij} })^*)_{kj}\enspace,
\quad\text{ for }n\geq n_{ij},
\end{equation}
where $q,q'$ are arbitrary numbers in $\fset{\sigma_{ij}}$
such that $q+q'\=mod n[\sigma_{ij}]$.
Furthermore, when $\sigma(A)$ is finite,
the integer $\sigma_{ij}$ can be chosen
so that it divides $\sigma(A)$.
\end{theorem}
The special case of this result when the state space $S$ is
finite was established in~\cite{cohen83}.
The book~\cite{bcoq} contains (Theorems~3.112 and~3.109) a partial account
of the results of that paper and is perhaps more easily available.

The integers $n_{ij}$ are sometimes called \new{coupling times}.
As  we discuss in Remark~\ref{rk-eff} below,
our proof of Theorem~\ref{cycl-th} 
yields explicit estimates of the coupling times, and of the integers
$\sigma_{ij}$. 
It also shows that every long
optimal path must pass through a critical node.
In fact, the following theorem
shows that optimal paths stay most of the time
in the critical graph.
\begin{theorem}[Turnpike theorem]\label{th-turnpike}
Make the same assumptions as in Theorem~\ref{cycl-th}.
Then, for all $i,j\in S$, there exists $m_{ij}\in\N$ having the 
following property:
any path from $i$ to $j$ having maximal weight amongst all paths 
with the same ends and length, has at most  $m_{ij}$ non-critical nodes.
\end{theorem}
The name ``turnpike theorem'' refers to a general class of 
results in dynamic programming, see for 
instance~\cite[Section 2.4]{maslovkololtsov95}.

To prove the cyclicity and turnpike theorems, 
we need a series of auxiliary results. 
\begin{lemma}\label{crit-inclu}
Assume that $\rho(A)=\unit$.
Let $i,j\in S$ and $\beta>\zero$, and write
$J:=\set{k\in S}{A^*_{ik} A^*_{kj}\geq \beta }$.
If $k\in J\cap N^r(A)$, then the recurrence class of $k$ is included
in $J$. In particular, if $k\in J\cap N^c(A)$, then the critical
class of $k$ is included in $J$.
\end{lemma}
\begin{proof}
Let $C$ be a recurrence class of $A$, and $k\in J\cap C$.
Then, $A^*_{k\ell} A^*_{\ell k}=\unit$ for all $\ell\in C$.
Hence, 
\[A^*_{i\ell} A^*_{\ell j}\geq A^*_{ik} A^*_{k\ell} A^*_{\ell k} A^*_{kj}
=A^*_{ik} A^*_{kj}\geq \beta\]
for all $\ell\in C$, which shows that $C\subset J$. 
This yields the  first assertion of the lemma.
Since critical classes are contained within recurrence classes, the
last assertion follows.
\end{proof}

\begin{corollary}
If $\tilde A$ has Property \PT,
then all the critical classes of $A$ are finite.
\end{corollary}
\begin{proof}
Let $C$ be a critical class of $A$.
Replace $A$ by $\tilde{A}$ and take
$i=j=k\in C$ and $\beta=\unit$ in Lemma~\ref{crit-inclu}.
Since $i\in J$, we get that $C\subset J$.
By Property \PT, $J$ is finite, hence $C$ is finite.
\end{proof}
We shall need the following standard result
of Perron-Frobenius theory.
\begin{lemma}\label{perron-graph}
Let $G$ be a (possibly infinite) strongly connected
graph and denote by $\gamma\in\N\setminus\{0\}$ its cyclicity.
Then, for all nodes $i$ and $j$ of $G$,
there exists $\nu_{ij}\in\N\setminus\{0\}$
such that:
\begin{itemize}
\item[-] all paths $p$ from $i$ to $j$ in $G$ 
satisfy $|p|\=mod\nu_{ij}[\gamma]$,
\item[-] for all $n\=mod\nu_{ij} [ \gamma]$ such that $n\geq \nu_{ij}$, 
there exists a path in $G$ from $i$ to $j$ of length $n$.
\end{itemize}
\end{lemma}
This lemma was proved in~\cite[Lemma~3.4.1 and~3.4.3]{brualdi}
when $G$ is finite, but the proof there also works
when $G$ is infinite. Similar results
were given in~\cite[Ch.~6,\S~3]{kemeny} in the context
of denumerable Markov chains. 

For any nodes $i$ and $j$ of a graph $G$, the least $\nu_{ij}$ satisfying 
the assertions of Lemma~\ref{perron-graph} will be denoted by $\nu_{ij}(G)$.
The problem of computing $\nu_{ij}(G)$ has been much studied
in Perron-Frobenius theory, where the maximum of $\nu_{ij}(G)$
over all $i,j$ is sometimes called the \new{exponent} of $G$.
It is known that when $\gamma=1$ and $G$ is finite,
$\nu_{ij}(G)\leq (|G|-1)^2+1$, where $|G|$
denotes the number of nodes of $G$; see~\cite[Th.~3.5.6]{brualdi}.
For an irreducible matrix $A$, we simply write $\nu_{ij}(A)$
instead of $\nu_{ij}(G(A))$.

For any matrices $A$ and integers $s,q$ such that $0\leq s< q$, 
we define $\tau_{ij}^{s,q}(A)$ to be the length of the
shortest path $p$ from $i$ to $j$ satisfying
$|p|\=mod s[q]$ and $|p|_A= (A^s (A^q)^*)_{ij}$.
By convention, $\tau_{ij}^{s,q}(A)=+\infty$
when no such path exists. The following
lemma gives a sufficient condition for the existence
of such a path.
\begin{lemma}\label{lem-newattained}
Assume that $A$ has Property \PT.
Then, for all $i,j\in S$ and $0\leq s<q$ such that
$(A^s(A^q)^*)_{ij}\neq\zero$, there is a path $p$ from 
$i$ to $j$ such that $|p|\=mod s[q]$ and $|p|_A=(A^s(A^q)^*)_{ij}$.
\end{lemma}
\begin{proof}
Consider the vector $v$ such that $v_k=A^s_{ik}(A^q)^*_{kj}$.
Since $v_k\leq A^*_{ik}A^*_{kj}$, 
Observations~\ref{obs-tight-leq} and~\ref{obs-2}
show that $v$ is tight. By Lemma~\ref{lem-supattained},
$(A^s(A^q)^*)_{ij}=A^s_{ik}(A^q)^*_{kj}$ for some $k\in S$.
By Proposition~\ref{lem-eq16}, $A^s_{ik}=|p|_A$ for some
path $p$ from $i$ to $k$ in the graph of $A$.
Now $(A^q)^*\leq A^*$, and so
$A^q$ also has Property \PT.
Using Proposition~\ref{lem-eq16} again, we get that
$(A^q)^*_{kj}=|p'|_{A^q}$ for some
some path $p'$ from $k$ to $j$
in the graph of $A^q$. Using Proposition~\ref{lem-eq16}
a third time, $|p'|_{A^q}=|p''|_{A}$ for some path $p''$
from $k$ to $j$ in the graph of $A$ whose length is a multiple of $q$.
The result follows.
\end{proof}
When $G(A)$ is finite,
we have the following explicit bound: 
$\tau_{ij}^{s,q}(A)\leq s+q (|G(A)|-1)$,
provided that $(A^s(A^q)^*)_{ij}\in\R$.

For each critical node $k$ of $A$, we denote by $\cycc{k}$ the cyclicity of the
strongly connected component of $k$ in $G^c(A)$. 
Observe that $\cycc{k}$ divides $\cycc{A}$ and, if $A$ is irreducible,
is a multiple of $\cyc{A}$.

The following technical lemma provides a lower
bound on the ultimate values of $A^n_{ij}$ as $n$ tends
to infinity. It will follow that any sufficiently long optimal
path must remain within the finite set $J$ defined in the lemma.
This will allow us to control the asymptotic behavior of $A^n_{ij}$
in essentially the same maner as in the finite dimensional case.
\begin{lemma}\label{lem-newnew}
Let $A$ be as in Theorem~\ref{cycl-th}, and assume furthermore that
$\rho(A)=\unit$.
Let $\gamma:=\cyc{A}$ be the cyclicity of $G(A)$.
Let $i,j\in S$ and $t\in\fset{\gamma}$ be such that
$t\equiv \nu_{ij}(A)\, [\gamma]$. 
Then, we can find at least one critical node $\ell$ such that
\[
\beta^\ell_{ij}:=
\min_{1\leq u\leq \cycc{\ell}/\gamma} (A^{t+u\gamma} (A^{\cycc{\ell}})^*)_{i\ell}
(A^{\cycc{\ell}})^*_{\ell j} \neq\zero \enspace .
\]
Choose such an $\ell$ arbitrarily, take any
number $\zero<\beta_{ij}\leq\beta_{ij}^\ell$,
define $J:=\set{k\in S}{A^*_{ik} A^*_{kj}\geq \beta_{ij}}$,
\(
\sigma_{ij}:= \lcm\set{\cycc{k}}{k\in N^c(A)\cap J}\),
and, for each $r\in\fset{\sigma_{ij}}$, define
\begin{equation}
\label{cycl-th1}
 Q^r_{ij}:=\bigoplus_{k\in N^c(A)}
\;
\bigoplus_{\scriptstyle s,s'\in \fset{{\sigma_{ij}}}\atop\scriptstyle s+s'\=mod r[\sigma_{ij}]}
(A^s (A^{\sigma_{ij}})^*)_{ik} (A^{s'} (A^{\sigma_{ij}} )^*)_{kj}
\enspace .
\end{equation}
Then,
\begin{align}
Q^{r}_{ij}\geq \beta_{ij} \text{ when }  r\=mod \nu_{ij}(A) [\gamma],
 \text{ and } 
Q^{r}_{ij}=\zero  \text{ otherwise.}
\label{cycl-min}
\end{align}
\end{lemma}
\begin{proof}
Take any critical node $k$. Since $A$ is irreducible, 
there exists a path $p_1$ from $k$ to $j$ in $G(A)$.
Let $m\in\fset{\cycc{k}}$ be such that $m+|p_1|\equiv 0\, [\cycc{k}]$.
Let $c$ be a critical circuit passing through $k$. Then, $|c|\geq 1$ and
by definition, $\cycc{k}$ divides $|c|$.
Moreover, there exists $\ell$ in the same strongly connected component as
$k$ in  $N^c(A)$,  such that $c$ can be written as the concatenation 
of a path $p_2$ from $k$ to $\ell$ of length $|c|-m$ and 
a path $p_3$ from $\ell$ to $k$ of length $m$.
Hence, $p_3 p_1$ is a path from $\ell$ to $j$
of length $m+|p_1|\equiv 0\, [\cycc{k}]$. Since
$\cycc{\ell}=\cycc{k}$, this implies that $(A^{\cycc{\ell}})^*_{\ell j}\neq \zero$.
By Lemma~\ref{perron-graph}, for all $n\geq \nu_{i\ell}(A)$ such that 
$n\equiv \nu_{i\ell}(A)\, [\gamma]$, $A^n_{i\ell}\neq \zero$.
Hence, if $t'\equiv \nu_{i\ell}(A)\, [\gamma]$,
$\min_{1\leq u\leq \cycc{\ell}/\gamma} 
(A^{t'+u\gamma} (A^{\cycc{\ell}})^*)_{i\ell}
(A^{\cycc{\ell}})^*_{\ell j}\neq \zero$
(since $\gamma$ divides $\cycc{\ell}$).
Then, necessarily $t'\equiv \nu_{ij}(A)\, [\gamma]$ (otherwise the above
term is $\zero$), which shows that $\beta^\ell_{ij}\neq \zero$.

Since $A$ has Property \PT, the set $J$ is finite.
By Lemma~\ref{crit-inclu}, $J$ contains any critical class
that it intersects. Since the critical node $\ell$ belongs
to $J$, we conclude that $\rho(A_{JJ})=\unit$ and 
$N^c(A)\cap J =N^c(A_{JJ})$. Therefore $\sigma_{ij}=\cycc{A_{JJ}}$.
Moreover, $\gamma$ divides $\sigma_{ij}$.

Fix $r\in\fset{\sigma_{ij}}$.
By definition, $Q^r_{ij}$ is the supremum of
all the weights $|p|_A$ of paths $p$ from $i$ to $j$ 
of length $|p|\equiv r\, [{\sigma_{ij}}]$ that pass through $N^c(A)$.
In particular, 
$Q^r_{ij}=\zero$ for all $r\not\=mod \nu_{ij}(A)[\gamma]$.
Assume now that $r\equiv \nu_{ij}(A)\, [\gamma]$.
Since $\ell\in N^c(A)\cap J$, $\sigma(\ell)$ divides $\sigma_{ij}$, 
and we can choose $u\in\{1,\ldots,\cycc{\ell}/\gamma\}$
such that $t+u\gamma \=mod r[\cycc{\ell}]$.
Therefore
$\beta_{ij}\leq (A^{t+u\gamma} (A^{\cycc{\ell}})^*)_{i\ell}
(A^{\cycc{\ell}})^*_{\ell j}$. 
By Lemma~\ref{lem-newattained},
we can find a path $p$ from $i$ to $\ell$ with $|p|\=mod t+u\gamma[\cycc{\ell}]$,
such that $(A^{t+u\gamma} (A^{\cycc{\ell}})^*)_{i\ell}=|p|_A$,
and we can find a path $p'$ from $\ell$ to $j$, with
$|p'|\=mod 0[\cycc{\ell}]$, such that $(A^{\cycc{\ell}})^*_{\ell j}=|p'|_A$.
Then, $\beta_{ij}\leq |p|_A|p'|_A$. 
By Lemma~\ref{perron-graph}, for any multiple $n$ of $\cycc{\ell}$ that is
large enough,
there exists a path $p''$ in $G^c(A)$ from $\ell$ to $\ell$ of length $n$.
The weight $|p''|_A$ of this path is $0$.
Since $\sigma(\ell)$ divides $\sigma_{ij}$,
we can choose $n$ such that $|p|+n+|p'|\=mod r[\sigma_{ij}]$.
Concatenating $p$, $p''$ and $p'$ we get a path
$p'''$ from $i$ to $j$ of length
$|p'''|\=mod r[\sigma_{ij}]$ passing through
$N^c(A)$. Thus $\beta_{ij}\leq |p|_A|p''|_A|p'|_A
= |p'''|_A\leq Q^r_{ij}$.
\end{proof}
\begin{lemma}\label{lem-cyc-new}
Let $A$, $\gamma$, $i$, $j$, $\beta_{ij}$, $J$, $\sigma_{ij}$, $r$, 
and $Q_{ij}^r$ be as in Lemma~\ref{lem-newnew}. 
Denote $J':=J\setminus N^c(A)$ and $B:=A_{J'J'}$, and define
\begin{align}
\mu_{ij}^r := \begin{cases}
0 &\text{if either } i\not\in J'\text{ or } j\not\in J'\enspace,\\
|J'| &\text{if } \rho(B)=\zero \text{ and } i,j\in J'\enspace,\\
\min\set{m\in\N}{
\rho(B)^{m}  \tilde{B}^*_{ij} \leq Q^r_{ij}}&\text{otherwise.}
\end{cases}
\label{e-muijr}
\end{align}
Denote by $N_{ij}^r$ the set of triples $(k,s,s')$ which attain the maximum
in~\eqref{cycl-th1}, and define
\begin{equation}\label{cycl-th3} 
\nu_{ij}^r:=\min_{(k,s,s')\in N_{ij}^r}
\tau_{ik}^{s,{\sigma_{ij}}}(A)+
{\sigma_{ij}} \nu_{kk}(G^c(A^{\sigma_{ij}}))+\tau_{kj}^{s',{\sigma_{ij}}}(A)
\enspace .
\end{equation}
Then
$A^n_{ij}=Q_{ij}^r$ for all
$n\geq n_{ij}^r$ such that $n\=mod r [\sigma_{ij}]$, where 
\(
n^r_{ij}:=\max(\mu_{ij}^r,\nu_{ij}^r)
\).
\end{lemma}
\begin{proof}
If $r\not\equiv \nu_{ij}(A)\, [\gamma]$, then, by Lemma~\ref{perron-graph}
and Equation~\eqref{cycl-min},
$Q^r_{ij}=A^n_{ij}=\zero$ for all $n\equiv r \, [{\sigma_{ij}} ]$.
So we assume that  $r\equiv \nu_{ij}(A)\, [\gamma]$.

Let us first show that $A_{ij}^n\leq Q^r_{ij}$ for all 
$n\=mod r[\sigma_{ij}]$ such that $n\geq \mu_{ij}^r$.
We saw in the proof of Lemma~\ref{lem-newnew} that
$\rho(A_{JJ})=\unit$ and $N^c(A)\cap J =N^c(A_{JJ})$.
Hence $\rho(B)<\unit$.
Let $n\equiv r\, [{\sigma_{ij}}]$. By Lemma~\ref{lem-eq16},
there exists a path $p$ from $i$ to $j$ of length $n$
such that $A^n_{ij}=|p|_A$.
If $p$ intersects $N^c(A)$, then, by definition of $Q^r_{ij}$,
$|p|_A\leq Q^r_{ij}$.
If $p$ intersects $S\setminus J$, then, by definition of $J$, 
$|p|_A< \beta_{ij}$. But $\beta_{ij} \leq Q^r_{ij}$,
by Lemma~\ref{lem-newnew},
and therefore, $|p|_A\leq Q_{ij}^r$. 
Otherwise, $p$ is included in $J'$, which implies that
$i,j\in J'$ and $|p|_A\leq B^n_{ij}$.
This shows that for all $n\=mod r[\sigma_{ij}]$, 
\begin{align*}
A^n_{ij}\leq Q^r_{ij}\oplus B^n_{ij} 
\text{ if } i,j\in J'\quad \text{ and }
A^n_{ij}\leq Q^r_{ij} \text{ otherwise.} 
\end{align*}
Since $Q^r_{ij}\neq\zero$, the integer $\mu_{ij}^r$ 
defined by~\eqref{e-muijr} exists. 
Applying again the observation of the proof 
of Theorem~\ref{cor-cycl1}, 
we get that $B^n=\zero$ for all $n\geq |J'|$, if $\rho(B)=\zero$,
and $B^n\leq \rho(B)^n\tilde{B}^*$ for all $n\geq 0$, if $\rho(B)>\zero$.
Hence, by definition of $\mu_{ij}^r$, we have
$A^n_{ij}\leq Q^r_{ij}$ for all $n\geq \mu_{ij}^r$
such that $n\=mod r[\sigma_{ij}]$.

We now show that the reverse inequality
holds for all
$n\=mod r[\sigma_{ij}]$ such that $n\geq \nu_{ij}^r$.
Since $(A^s (A^{\sigma_{ij}})^*)_{ik} (A^{s'} (A^{\sigma_{ij}} )^*)_{kj}
\leq A^*_{ik} A^*_{kj}$, and since $Q^r_{ij}\geq \beta_{ij}\neq\zero$, the set
$N_{ij}^{r}$ is included in $J\times\fset{\sigma_{ij}}^2 $,
and Lemma~\ref{lem-supattained}
shows that $N_{ij}^r\neq\emptyset$. 
This implies that the integer $\nu_{ij}^{r}$ defined by~\eqref{cycl-th3}
exists.
Moreover, $\nu_{ij}^{r}\equiv r\, [{\sigma_{ij}}]$.
Let $(k,s,s')\in N_{ij}^{r}$ be a triple 
attaining the minimum in~\eqref{cycl-th3}.
By Lemma~\ref{lem-newattained},
we can find a path $p$ from $i$ to $k$ of length $\tau_{ik}^{s,\sigma_{ij}}(A)$
such that $(A^{s} (A^{\sigma_{ij}})^*)_{ik}=|p|_A$,
and a path $p'$ from $k$ to $j$ of length
$\tau_{k j}^{s',\sigma_{ij}}(A)$ such that
$(A^{s'}(A^{\sigma_{ij}})^*)_{k j}=|p'|_A$.
Since $\cycc{k}$ divides $\sigma_{ij}$, then 
$\rho(A^{\sigma_{ij}})=\unit$, and since to every
path of length $n\sigma_{ij}$ in $G^c(A)$ corresponds
a path of length $n$ in $G^c(A^{\sigma_{ij}})$,
it follows from Lemma~\ref{perron-graph} that
the cyclicity of the strongly connected component
of $k$ in $G^c(A^{\sigma_{ij}})$ is equal to $1$.
Applying Lemma~\ref{perron-graph} again,
we get that for each $n\geq \nu_{kk}(G^c(A^{\sigma_{ij}}))$,
there is a path $p''$ in $G^c(A^{\sigma_{ij}})$ from $k$ to $k$ of length
$n$. To this path corresponds a path $p'''$ in $G^c(A)$ 
of length $n\sigma_{ij}$, with the same ends.
We have $|p'''|_A=\unit$. Concatenating $p$, $p'''$ and $p'$ we get 
a path from $i$ to $j$ of length
$|p|+n\sigma_{ij}+|p'|\=mod r[\sigma_{ij}]$
such that $Q^r_{ij}= |p|_A |p'''|_A|p'|_A
\leq A^{|p|+n\sigma_{ij}+|p'|}_{ij}$. 
We deduce that for all $n\geq \nu_{ij}^r$
such that $n\=mod r[\sigma_{ij}]$, 
$Q^r_{ij}\leq A^n_{ij}$.
Then, for all $n\geq \max(\nu_{ij}^{r},\mu_{ij}^r)$
such that  $n\equiv r \, [{\sigma_{ij}}]$, $Q^r_{ij}= A^n_{ij}$. 
\end{proof}
\begin{lemma}\label{lem-new}
Let $A$, $i$, $j$, $\sigma_{ij}$,
$r$, and $Q_{ij}^r$ be as in Lemma~\ref{lem-newnew}.
For all $q,q'\in\fset{\sigma_{ij}}$ such that 
$q+q'\=mod r[\sigma_{ij}]$, we have
\begin{align}
Q^r_{ij}= 
\bigoplus_{k\in N^c(A)}
(A^q (A^{\sigma_{ij}})^*)_{ik} (A^{q'} (A^{\sigma_{ij}} )^*)_{kj}
\enspace .\label{e-ambig}
\end{align}
\end{lemma}
\begin{proof}
We denote by $P_{ij}$ the right hand side of~\eqref{e-ambig},
and write $\gamma$ instead of $\gamma(A)$, as before.
Trivially, $P_{ij}\leq Q^r_{ij}$,
hence, by~\eqref{cycl-min}, $P_{ij}=Q^r_{ij}=\zero$ 
when $r\not\=mod \nu_{ij}(A) [\gamma]$.
So we assume that $r\=mod \nu_{ij}(A)[\gamma]$.
Choose $(k,s,s')$ in the set $N_{ij}^r$
defined in the proof of Lemma~\ref{lem-cyc-new}.
Take also $p,p',p'''$ as in the same proof,
with $|p'''|=\sigma_{ij}\nu_{kk}(G^c(A^{\sigma_{ij}}))$.
We have $Q_{ij}^r=|p|_A|p'''|_A|p'|_A$.
Moreover, for all $m\in\{0,\ldots,|p'''|\}$, and
in particular, for $0\leq m< \sigma_{ij}$,
there exists $\ell$ in $N^c(A)$
such that $p'''$ can be written as the concatenation 
of a path $p_1$ from $k$ to $\ell$ of length $m$ and 
a path $p_2$ from $\ell$ to $k$ of length $|p'''|-m$.
The concatenated path $p p_1$ goes from $i$ to $\ell$
and has length $|p|+m\equiv s+m\, [\sigma_{ij}]$.
The concatenated path $p_2p'$ goes from $\ell$ to $j$
and has length $|p'''|-m+|p'|\=mod -m+s'[\sigma_{ij}]$.
Choose $m$ such that $s+m\=mod q[\sigma_{ij}]$,
so that $-m+s'\=mod q'[\sigma_{ij}]$. Then,
$|pp_1|_A\leq (A^q (A^{\sigma_{ij}})^*)_{i\ell}$ and
$|p_2p'|_A\leq (A^{q'} (A^{\sigma_{ij}} )^*)_{\ell j}$.
This shows that $Q_{ij}^r=|p|_A|p'''|_A|p'|_A
=|pp_1|_A |p_2p'|_A \leq P_{ij}$.
\end{proof}
\begin{proof}[Proof of Theorem~\ref{cycl-th}]\sloppy
We can assume without loss of generality that $\rho(A)=\unit$, 
so that $\tilde{A}=A$.
We obtain Equation~\eqref{cycl1} from Lemmas~\ref{lem-cyc-new}
and~\ref{lem-new} by taking  
$\sigma_{ij}$ as in Lemma~\ref{lem-newnew} and  
\begin{align*}
n_{ij}:=\max_{0\leq r\leq \sigma_{ij}-1}n_{ij}^r\enspace.
\end{align*}
Equation~\eqref{cycl1p3} follows immediately.
\end{proof}
\begin{proof}[Proof of Theorem~\ref{th-turnpike}]
We assume again, without loss of generality, that $\rho(A)=\unit$,
so that $\tilde{A}=A$. Let $n\geq 1$ and let $p_1$ be a path of length
$n$ from $i$ to $j$ that has maximal weight amongst all such paths,
 that is, $|p_1|_A=A^n_{ij}>\zero$.
Let $\sigma_{ij}$ and $n_{ij}$ be defined as in Theorem~\ref{cycl-th}.
Let $\beta_{ij}$ and $J$ be defined as in Lemma~\ref{lem-newnew},
and, for $r\in \fset{\sigma_{ij}}$,
with $r\=mod n[\sigma_{ij}]$,
let $Q^r_{ij}$ be defined
as in Lemma~\ref{lem-newnew}.
By Theorem~\ref{cycl-th}, Lemma~\ref{lem-new}
and Lemma~\ref{lem-newnew},
we have $A^n_{ij}=Q^r_{ij}\geq \beta_{ij}$ if $n\geq n_{ij}$.
So, all the nodes of $p_1$ must belong to $J$. 
For each path $p$ in $J$, let $\swei{p}$ denote the number
of arcs of $p$ whose initial node is non-critical.
In particular, when $p$ is a circuit, 
$\swei{p}=0$ if and only if all the nodes of $p$
are critical.
Let $\vartheta$ denote the maximal value of 
the mean weight $|c|_A^{1/{\swei{c}}}$
over all elementary circuits $c$ that stay within $J$ and contain
at least one non-critical node. Since $J$ is finite,
$\vartheta<\unit$. We next show that, when $n\geq n_{ij}$,
the number of non-critical nodes of $p_1$ can be bounded
by a constant $\kappa_{ij}^r$.
Write $p_1$ as a disjoint union
of an elementary
path $p_0$ from $i$ to $j$ and of elementary circuits $c_1,\ldots, c_s$. 
First consider the case where $\vartheta=\zero$. Then, all the circuits
that stay in $J$ are critical, and so the number
of non-critical nodes of $p_1$ is at most equal to the number
of nodes of $p_0$, which is at most equal to
the number of nodes of $J$, since $p_0$ is elementary. 
Thus, in this case,
we can take $\kappa_{ij}^r:=|J|$.
Suppose now that $\vartheta\neq \zero$.
Since the weight of any circuit contained in $J$
is at most $\vartheta^{\swei{c}}$,
we have 
\[ A^n_{ij}=|p_1|_A\leq |p_0|_A\vartheta^{\swei{c_1}+\cdots +\swei{c_s}}
=|p_0|\vartheta^{\swei{p_1}-\swei{p_0}}\leq C_{ij}\vartheta^{\swei{p_1}}
\enspace,\]
where $C_{ij}$ is equal to the maximum of 
the quantity $|p'_0|_A\vartheta^{-\swei{p'_0}}$ over
all the elementary paths $p'_0$ from $i$ to $j$ that stay
in $J$. Since the number of non-critical nodes of $p_1$
is at most $\swei{p_1}+1$,
we can take $\kappa_{ij}^r$ to be the least integer such that
$C_{ij}\vartheta^{\kappa_{ij}^r}< Q^r_{ij}$.
Then, the number of non-critical nodes of $p_1$ is bounded
above by
$m_{ij}:=\max_{r\in \fset{\sigma_{ij}}} \max(n_{ij},\kappa^r_{ij})$.
\end{proof}
\begin{remark}
The quantities $\vartheta$ and $C_{ij}$ which appears in the proof of Theorem~\ref{th-turnpike} can be computed by max-plus Schur complement formulas, see~\cite[Section~4]{abg04a}.
\end{remark}
\begin{remark}\label{rk-eff}
The proof of Theorem~\ref{cycl-th} yields an explicit estimate of the
integers $\sigma_{ij}$ and the coupling times
$n_{ij}$, appearing in Theorem~\ref{cycl-th}.
The bound on the coupling time involves an essentially ``arithmetical''
term, $\nu_{ij}^r$ (depending on the constants $\nu_{kk}$),
and the maximal circuit mean $\rho(B)$ of a certain submatrix $B$
of $A$ (see Lemma~\ref{lem-cyc-new}). In the case when $S$ is finite, 
the problem of estimating $n_{ij}$ has received attention.
In this special case, the bound of the present
paper improves the bound of~\cite{BG01}, which uses related ideas.
A bound of a different nature has appeared in~\cite{hartmann}.
\end{remark}
\begin{remark}
The bound of $n_{ij}$ depends
on the choice of the critical node $\ell$ in Lemma~\ref{lem-newnew}.
Choosing the critical node $\ell$ which maximizes the quantity 
$\beta_{ij}^\ell$
defined in Lemma~\ref{lem-newnew} yields the smallest set $J$, and therefore, 
the best constants $\sigma_{ij}$ and $n_{ij}$.
When $S$ is finite, the technicalities
of Lemma~\ref{lem-newnew} might be dispensed with,
at the price of a coarser bound. We may just take $\beta_{ij}=\zero$ and
$J:=S$ in Lemma~\ref{lem-newnew} and Lemma~\ref{lem-cyc-new}
to compute $\sigma_{ij}$ and $n_{ij}$.
\end{remark}
\begin{remark}
Example~\ref{ex-rec-noncrit2} shows that Property \PT\ 
cannot be dispensed with in Theorem~\ref{cor-cycl1}.
In this example, $\rho(A)=\unit$ and $N^c(A)=N^r(A)=\emptyset$,
but $A^n_{00}$ does not tend to $\zero$ as $n$ tends to
infinity. 
\end{remark}
Without Property \PT, one would not expect convergence 
in finite time to a periodic regime as in Theorem~\ref{cycl-th}.
The following counter-example
shows that one might not even have asymptotic convergence.
\begin{example}
We will give an irreducible matrix $A$ with the following
features: $A$ does not have Property \PT,
$N^c(A)=N^r(A)\neq\emptyset$, $\rho(A)=\unit$,
the sequence $(A^n_{11})_{n\geq 2}$ is bounded
but there is no positive integer $\sigma$ such that
$A^{n\sigma}_{11}$ converges as $n$ tends to infinity.
Let $S:=\N$, and let $\alpha_2,\alpha_3,\ldots$ denote any sequence
of negative numbers.
We set $A_{00}:=0$, $A_{01}=A_{10}:=-1$, $A_{i,i+1}:=0$
and $A_{i+1,1}:=\alpha_{i+1}$ for $i\geq 1$.
All other arcs are given weight $-\infty$.
The graph of $A$ is 
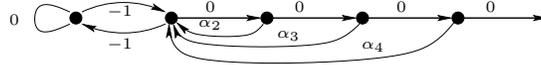
\begin{center}
\input nocritical5
\end{center}
If $\alpha_n\geq -2$
and $\alpha_{n+p}\geq \alpha_n +\alpha_p$
for all $n,p\geq 2$,
then $A^n_{11}=\alpha_n$ for $n\geq 2$. 
Choosing for $\alpha_2,\alpha_3,\ldots$ any sequence
taking only the values $-1$ and $-2$ and having
arbitrarily many consecutive occurences of each,
we get the announced property for the sequence $A^n_{11}$.
\end{example}
\section{Representation of Max-Plus Eigenvectors and 
Max-Plus Martin Boundary}\label{s-martin}
In this section, we present, without proof, 
some of the results of~\cite{AGW-m},
and show how they relate to some of the examples we have encountered.
The focus of attention in that paper was the $\unit$-eigenspace.
By analogy with potential theory,
the elements of this eigenspace were called \new{harmonic vectors} and
elements of the $\unit$-super-eigenspace were called
\new{super-harmonic vectors}.
Here, we restate the results in terms of a general $\lambda$-eigenspace,
a trivial change since this eigenspace is exactly the set of harmonic vectors
with respect to $\Al:= \lambda^{-1} A$.

We shall make the following assumption.
\begin{assumption}
There exists a row $\lambda$-super-eigenvector with full support,
that is a row vector $\pi\in \R^S$ such that $\lambda \pi\geq \pi A$.
\end{assumption}
This assumption implies, in particular, that $\rho(A)\leq \lambda$
(see Lemma~\ref{l1}), that $\pi=\pi (\Al )^*$, and that 
$(\Al )^*_{ij}\in\rmax$ for all $i,j\in S$.

We shall look for eigenvectors $u$ that are \new{$\pi$-integrable},
meaning that $\pi u<+\infty$.
We denote by $\El$ the boundedly complete subsemimodule
of $\rmax^S$ consisting of those $\pi$-integrable vectors $u$
such that $Au=\lambda u$. 

It is often possible to choose
$\pi:=(\Al )^*_{b\cdot}$ for some $b\in S$,
for instance when $A$ is irreducible.
With this choice, every $\lambda$-eigenvector $u\in\rmax^S$ 
is automatically $\pi$-integrable.

We define the \new{Martin kernel} $K$ with respect to $\pi$ and $\lambda$:
\begin{align*}
K_{ij} := (\Al )^*_{ij}(\pi_{j})^{-1} \quad \mbox{for all $i,j \in S$} \enspace .
\end{align*}
Since $\pi_{i}(\Al )^*_{ij}\leq (\pi (\Al )^*)_{j}=\pi_{j}$, we have 
\begin{align*}
K_{ij}\leq (\pi_{i})^{-1} \quad \mbox{for all $i,j \in S$} \enspace.
\end{align*}
This shows that the columns $K_{\cdot j}$ are bounded above
independently of $j$. By Tychonoff's theorem, the set of columns
$\sK= \set{K_{\cdot j}}{j\in S}$ is relatively compact 
in the product topology of $\rmax^S$. The \new{Martin space}
$\sM$ is defined to be the closure of $\sK$ in this space.
We call $\sB:=\sM\setminus \sK$ the \new{Martin boundary}.

Let $u\in \rmax^S$ be a $\pi$-integrable vector.
We define the map $\mu_u:\sM\to \rmax$ by
\begin{align*}
\mu_u(w) := \limsup_{K_{\cdot j} \to w} \pi_j {u}_j
 := \inf_{W\ni w} \sup_{K_{\cdot j} \in W} \pi_j {u}_j 
\quad\text{for $w\in \sM$} \enspace,
\end{align*}
where the infimum is taken over all open neighborhoods $W$ of $w$ in $\sM$. 
The map $\mu_u$ is automatically upper semicontinuous
and bounded above by $\pi u <+\infty$.

We wish to define a particular subset of the Martin space,
called the \new{minimal Martin space}. To do this, we introduce a
kernel $H^{\pl}$ over $\sM$ which extends, in some sense, the $(\Al )^+$ 
matrix:
\begin{align*}
H^{\pl}(w',w) & := 
\limsup_{K_{\cdot i}\to w'} \, \liminf_{K_{\cdot j}\to w} \pi_i 
(\Al )^+_{ij} (\pi_j)^{-1} \enspace .
\end{align*}
Note that
\begin{align*} 
H^{\pl}(w',w)\leq \unit \quad \mbox{for all $w,w'\in \sM$} \enspace.
\end{align*}
When $i,j\in S$,
\begin{align*}
H^{\pl}(K_{\cdot i},K_{\cdot j})&=\pi_i (\Al )^+_{ij} 
(\pi_j)^{-1}
\enspace .
\end{align*}
We now define the \new{minimal Martin space} to be
\begin{align*}
 \sMin:= \set{w\in \sM}{H^{\pl}(w,w)=\unit}\enspace .
\end{align*}
\begin{theorem}[Poisson-Martin representation, {\cite{AGW-m}}]
\label{poisson-martin}
Any element $u\in\El$ can be written as
\begin{align}
u=\bigoplus_{w\in \sMin} \nu(w) w \enspace ,
\label{Hequal1}
\end{align}
with $\nu:\sMin\to\rmax$, and necessarily,
\begin{align}
\sup_{w\in \sMin} \nu(w) <+\infty \enspace .
\label{Hequal1b}
\end{align}
Conversely, any $\nu:\sMin\to\rmax$
satisfying~\eqref{Hequal1b} 
defines by~\eqref{Hequal1} an element $u$ of $\El$.
Moreover, given $u\in \El$,
$\mu_u$ is the maximal $\nu$ satisfying~\eqref{Hequal1}.
\end{theorem}

This theorem shows, in particular, that $\El=\{\zero\}$
if and only if $\sMin$ is empty.

We say that a vector $u\in \rmax^S$ is \new{normalized} if $\pi u=\unit$.
If a subsemimodule of $\rmax^S$ contains only $\pi$-integrable
vectors, then its extremal generators are
exactly those vectors of the form $\alpha \xi$, with $\alpha\in \R$ and
$\xi$ a normalized extremal generator.

\begin{theorem}[{\cite{AGW-m}}]
The normalized extremal generators of $\El$ are precisely the elements
of $\sMin$.
\end{theorem}
\begin{remark}\label{rk-finrec}
Suppose $\lambda=\rho(A)$, there are only finitely many recurrence classes,
and all but finitely many nodes are recurrent.
In this case, $\sK$ is a finite set, and so $\sM=\sK$, the boundary $\sB$
is empty, and $\sMin$ is the set of columns
$K_{\cdot j}$, with $j$ recurrent. Then, the
representation theorem (Theorem~\ref{poisson-martin})
shows that any $\rho(A)$-eigenvector is a finite linear combination of 
the recurrent columns of $\tilde{A}^*$, just as in the finite dimensional
case. We saw an example 
of this situation in Examples~\ref{ex-rec-noncrit}
and~\ref{ex-rec-noncrit3}. There $S=\N$ was the only recurrence class,
and so any $\unit$-eigenvector had to be a multiple of 
$\tilde{A}^*_{\cdot 0}\equiv \unit$, and hence constant.
\end{remark}

In~\cite{AGW-m}, we also prove a representation theorem for 
super-eigenvectors similar to Theorem~\ref{poisson-martin}, with
$\sMin$ replaced by $\sMin\cup\sK$.
Moreover, we characterize $\sMin\cup\sK$
as the set of normalized extremal generators of the set of $\pi$-integrable
vectors satisfying $Au\leq \lambda u$.

The following result gives a  condition which guarantees 
the existence of eigenvectors.

\begin{prop}[{\cite{AGW-m}}]\label{prop-cs}
Assume that $S$ is infinite, that 
the vector $\pi^{-1}:=((\pi_i)^{-1})_{i\in S}$
is $A$-tight and that $\zero\not\in\sM$.
Then, $\sMin$ is non-empty.
\end{prop}

\begin{corollary}[{\cite{AGW-m}}]\label{cor-a-irred}
Assume that $S$ is infinite, and that
$A$ is irreducible and right locally finite. Then
the spectrum of $A$ is $[\rho(A),+\infty)$.
\end{corollary}
\begin{example}
The matrix $A$ of Example~\ref{ex-birth}, with $p,q\neq \zero$, satisfies
the assumptions of Corollary~\ref{cor-a-irred},
and so its spectrum is $[\rho(A),+\infty)$,
as we have already noted.
Let us compute the $\lambda$-eigenspace of $A$ for 
$\lambda\geq \rho(A)=\sqrt{pq}$.
We take $\pi:=(\Al)^*_{0\cdot}$.
Clearly,
$\pi_j=({p}/{\lambda})^j$ for all $j\in\N$.
Also, $K_{ij}= ({\lambda}/{p})^i$ when $i\leq j$ and
$K_{ij}= ({\lambda}/{p})^i({pq}/{\lambda^2})^{i-j}$ otherwise.
Hence, $\sMin=\sB=\{u\}$, where the vector $u\in\R^\N$ is given by
$u_i:= ({\lambda}/{p})^i$ for all $i\in\N$.
So, Theorem~\ref{poisson-martin} shows that every $\lambda$--eigenvector 
is a multiple of $u$.
\end{example}
\begin{example} The matrix $A$ of Example~\ref{ex-triang}
is irreducible and  left locally finite,
so Corollary~\ref{cor-a-irred} shows
that the transpose of $A$ has an eigenvector
for all $\lambda \geq \rho(A) \enspace$.
However, $A$ has no eigenvectors, which shows that the assumption
that $A$ is right locally finite is needed in Corollary~\ref{cor-a-irred}.
One can also prove that $A$ has no eigenvectors,
using Theorem~\ref{poisson-martin}.
Indeed, let $\lambda\geq \unit$ and consider $\pi:=(\Al)^*_{0\cdot}$.
We have $\pi_0=0$ and $\pi_i=\beta -\lambda$ when $i\geq 1$.
Therefore, $\pi^{-1}$ is not $A$-tight, and so does not satisfy the conditions
of Proposition~\ref{prop-cs}. We have $K_{ij}=0$ when $i<j$, 
$K_{ij}= (\beta-\lambda)(i-j-1)$ when $i\geq j\geq 1$ and
$K_{ij}= (\beta-\lambda)i$ otherwise.
Hence, $\sB=\{\unit\}$, where $\unit$ is the unit vector.
Since $A$ has no recurrent nodes, $\sMin\subset \sB$, and since $\unit$ is not
an eigenvector of $A$, we deduce that $\sMin$ is empty.
Theorem~\ref{poisson-martin} then shows that, for all $\lambda\geq \rho(A)$,
there is no $\lambda$-eigenvector.
\end{example}
\begin{remark}
In general,
the $\rho(A)$-eigenspace $\sE$ of $A$ may or may not be spanned by the 
critical columns of $\tilde{A}^*$.
Let us first consider examples of matrices $A$ such that $\tilde{A}$ 
does not have Property \PT. 
It follows from Remark~\ref{rk-finrec}
that for the matrix $A$ of Example~\ref{ex-rec-noncrit}, $\sE$
is the set spanned by the unit vector
(the vector identically equal to $\unit$), whereas
the graph of $A$ has no critical nodes.
When $A$ is the matrix of Example~\ref{ex-rec-noncrit3}, $\sE$ is again
spanned by the unit vector, 
all the columns of $\tilde{A}^*$ are critical, 
and coincide with the unit vector,
and critical classes are singletons,
so that picking one column $\tilde{A}^*_{\cdot j}$ per critical
class of $A$, as in Theorem~\ref{generate-astar},
yields a non minimal spanning family of $\sE$.

Let us now consider examples of matrices $A$ such that $\tilde{A}$ 
has Property \PT.  Even in this case, 
there may exist $\rho(A)$-eigenvectors that do not have the 
representation~\eqref{e-rep}.
For example, the matrix in Example~\ref{ex-assum-tight2} has
no critical classes, whereas one can show that it has $\rho(A)$-eigenvectors.
Indeed, take $\pi=A^*_{0\cdot}\equiv \unit$.
Observing that $K_{ij}=A^*_{ij}\to \unit$ when $j$ goes to
infinity, we see that the boundary $\sB$ consists of just the unit vector
$\unit$. Since $A$ has no recurrence classes,
$\sMin\subset\sB$, and since $\unit$ is an eigenvector,
the minimal boundary is given by $\sMin=\sB=\{\unit\}$.
Thus by Theorem~\ref{poisson-martin}, every $\rho(A)$-eigenvector is constant.

To take another example, 
the matrix in Example~\ref{ex-assum-tight1} has a single
critical class $\{0\}$. The associated eigenvector 
$A^*_{\cdot 0}=(-i)_{i\in\N}$ is $A^*$-tight.
But, taking $\pi=A^*_{0\cdot}\equiv \unit$ as before,
we again get that $\sB=\{\unit\}$. Then, $\sMin=\{A^*_{\cdot 0},\unit\}$
and by Theorem~\ref{poisson-martin}, any $\unit$-eigenvector is a max-plus 
linear combination of $A^*_{\cdot 0}$ and $\unit$.
\end{remark}

\def\cprime{$'$}

\end{document}